%% file: m7-6.tex
\documentclass{gtart_h}

\input gtmonout

\volumenumber{7}
\volumename{Proceedings of the Casson Fest}
\volumeyear{2004}
\papernumber{6}
\pagenumbers{145}{180}
\received{14 August 2003}
\revised{9 January 2004}
\accepted{28 January 2004}
\published{18 September 2004}

\usepackage{amsmath,amssymb,epsf}
\usepackage[all]{xy}

\newtheorem{thm}{Theorem}[section]
\newtheorem*{thm*}{Theorem}
\newtheorem{Definition}[thm]{Definition}
\newtheorem{Example}[thm]{Example}
\newtheorem{Remark}[thm]{Remark}
\newtheorem{Proposition}[thm]{Proposition}

\newtheorem{Theorem}[thm]{Theorem}
\newtheorem{Lemma}[thm]{Lemma}
\newtheorem*{att}{Algebraic Transversality Theorem}
\newtheorem*{ctt}{Combinatorial Transversality Theorem}

\newcommand{\RR}{\mathbb R}
\newcommand{\Z}{\mathbb Z}

\begin{document}

\title{Algebraic and combinatorial\\codimension--1 transversality}
\asciititle{Algebraic and combinatorial codimension-1 transversality}

\author{Andrew Ranicki}

\address{School of Mathematics, University of Edinburgh\\King's Buildings,
Edinburgh EH9 3JZ, Scotland, UK}

\email{a.ranicki@ed.ac.uk}


\begin{abstract}
The Waldhausen construction of Mayer--Vietoris splittings of chain
complexes over an injective generalized free product of group
rings is extended to a combinatorial construction of Seifert--van
Kampen splittings of CW complexes with fundamental group an
injective generalized free product.
\end{abstract}
\asciiabstract{%
The Waldhausen construction of Mayer-Vietoris splittings of chain
complexes over an injective generalized free product of group rings is
extended to a combinatorial construction of Seifert-van Kampen
splittings of CW complexes with fundamental group an injective
generalized free product.}
\maketitlepage

\hfill{\small\it Dedicated to Andrew Casson}\rightskip25pt

\section*{Introduction}\rightskip0pt

The close relationship between the topological properties of
codimension--1 submanifolds and the algebraic properties of
groups with a generalized free product structure
first became apparent with the Seifert--van Kampen Theorem
on the fundamental group of a union, the work of Kneser
on 3--dimensional manifolds with fundamental group a free
product, and the topological proof of Grushko's theorem by
Stallings.

This paper describes two abstractions of the geometric codimension--1
transversality properties of manifolds (in all dimensions):
\begin{enumerate}
\item the algebraic transversality construction of Mayer--Vietoris
splittings of chain complexes of free modules
over the group ring of an injective generalized free product,
\item the combinatorial transversality construction of Seifert--van Kampen
splittings of CW complexes with fundamental group an injective generalized
free product.
\end{enumerate}

By definition, a group $G$ is a {\it generalized free product} if
it has one of the following structures:
\begin{itemize}
\item[{\rm (A)}] $G=G_1*_HG_2$ is the amalgamated free product
determined by group morphisms $i_1\co H \to G_1$, $i_2\co H \to G_2$, so
that there is defined a pushout square of groups
$$\xymatrix{H \ar[r]^-{i_1} \ar[d]_-{i_2} & G_1 \ar[d]^-{j_1}\\
G_2 \ar[r]^-{j_2} & G}$$
The amalgamated free product is {\it injective} if $i_1,i_2$ are
injective, in which case so are $j_1,j_2$, with
$$G_1 \cap G_2=H \subseteq G.$$
\item[] An injective amalgamated free product is {\it nontrivial} if
the morphisms $i_1\co H \to G_1$, $i_2\co H \to G_2$ are not isomorphisms,
in which case the group $G$ is infinite, and $G_1,G_2,H$ are subgroups
of infinite index in $G$.
\item[] The amalgamated free product is {\it finitely presented} if
the groups $G_1,G_2,H$ are finitely presented, in which case so is $G$.
(If $G$ is finitely presented, it does not follow that $G_1,G_2,H$
need be finitely presented).
\item[{\rm (B)}] $G=G_1*_H\{t\}$ is the HNN extension determined by
group morphisms $i_1,i_2\co H \to G_1$
$$H  \raise5pt\hbox{\xymatrix@R-25pt@C-5pt{\ar[r]^-{i_1}&\\ \ar[r]_{i_2}&}}
\xymatrix{ G_1 \ar[r]^-{j_1} & G}$$
with $t \in G$ such that
$$j_1i_1(h)t=tj_1i_2(h) \in G~~(h \in H).$$
The HNN extension is {\it injective} if $i_1,i_2$ are
injective, in which case so is $j_1$, with
$$G_1 \cap tG_1t^{-1}=i_1(H)=ti_2(H)t^{-1} \subseteq G$$
and $G$ is an infinite group with the subgroups
$G_1,H$ of infinite index in $G=G_1*_H\{t\}$.

\item[] The HNN extension is {\it finitely presented} if
the groups $G_1,H$ are finitely presented, in which case so is $G$.
(If $G$ is finitely presented, it does not follow that $G_1,H$
need be finitely presented).
\end{itemize}

A subgroup $H \subseteq G$ is {\it 2--sided} if $G$ is either an injective
amalgamated free product $G=G_1*_HG_2$ or an injective HNN extension
$G=G_1*_H\{t\}$.  (See Stallings \cite{stallings2} and Hausmann
\cite{hausmann} for the characterization of 2--sided subgroups in terms
of bipolar structures.)

A CW pair $(X,Y \subset X)$ is {\it 2--sided} if $Y$ has an open
neighbourhood $Y \times \RR \subset X$.  The pair is {\it connected} if
$X$ and $Y$ are connected.  By the Seifert--van
Kampen Theorem $\pi_1(X)$ is a generalized free product:
\begin{itemize}
\item[{\rm (A)}] if $Y$ separates $X$ then $X-Y$ has two components, and
$$X=X_1\cup_YX_2$$
for connected $X_1,X_2 \subset X$ with
$$\pi_1(X)=\pi_1(X_1)*_{\pi_1(Y)}\pi_1(X_2)$$
the amalgamated free product determined by the morphisms
$i_1\co \pi_1(Y) \to \pi_1(X_1)$, $i_2\co \pi_1(Y) \to \pi_1(X_2)$ induced by the
inclusions $i_1\co Y \to X_1$, $i_2\co Y \to X_2$.
$$
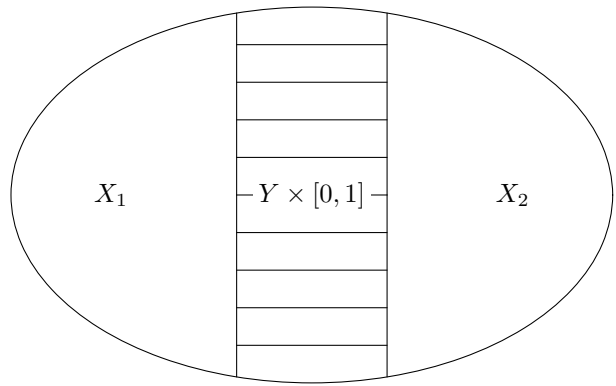
\epsfbox{a.ps}$$
\item[{\rm (B)}] if $Y$ does not separate $X$ then
$X-Y$ is connected and
$$X=X_1 \cup_{Y \times \{0,1\}} Y\times [0,1]$$
for connected $X_1 \subset X$, with
$$\pi_1(X)=\pi_1(X_1)*_{\pi_1(Y)}\{t\}$$
the HNN extension determined by the morphisms
$i_1,i_2\co \pi_1(Y) \to \pi_1(X_1)$
induced by the inclusions $i_1,i_2\co Y \to X_1$.
$$
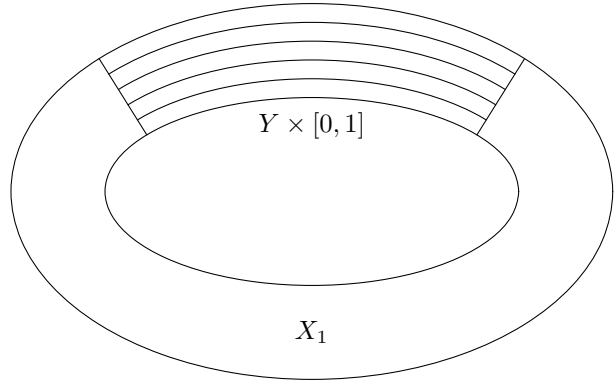
\epsfbox{b.ps}$$
\end{itemize}

The generalized free product is injective if and only if the
morphism $\pi_1(Y) \to \pi_1(X)$ is injective, in which case $\pi_1(Y)$
is a 2--sided subgroup of $\pi_1(X)$. In section \ref{tree} the
Seifert--van Kampen Theorem in the injective case will be deduced from
the Bass--Serre characterization of an injective $\begin{cases}
\hbox{\rm amalgamated free product} \\
\hbox{\rm HNN extension} \end{cases}$ structure on a group $G$ as
an action of $G$ on a tree $T$ with quotient
$$T/G=\begin{cases}
[0,1]\\
S^1~~~~.
\end{cases}$$

A codimension--1 submanifold $N^{n-1} \subset M^n$ is 2--sided
if the normal bundle is trivial, in which case $(M,N)$ is a
2--sided CW pair.

For a 2--sided CW pair $(X,Y)$ every map $f\co M \to X$ from an
$n$--dimensional manifold $M$ is homotopic to a map (also denoted by
$f$) which is transverse at $Y \subset X$, with
$$N^{n-1}=f^{-1}(Y) \subset M^n$$
a 2--sided codimension--1 submanifold, by the Sard--Thom theorem.

By definition, a {\it Seifert--van Kampen splitting} of a connected
CW complex $W$ with $\pi_1(W)=G=\begin{cases} G_1*_HG_2 \\ G_1*_H\{t\} \end{cases}$
an injective generalized free product is a connected 2--sided CW
pair $(X,Y)$ with a homotopy equivalence $X \to W$ such that
$${\rm im}(\pi_1(Y) \to \pi_1(X))=H \subseteq \pi_1(X)=\pi_1(W)=G.$$
The splitting is {\it injective} if $\pi_1(Y) \to \pi_1(X)$ is
injective, in which case
$$X=\begin{cases} X_1 \cup_Y X_2 \\
X_1 \cup_{Y \times \{0,1\}}Y \times [0,1]
\end{cases}$$
with
$$\begin{cases}
\pi_1(X_1)=G_1,\pi_1(X_2)=G_2\\
\pi_1(X_1)=G_1
\end{cases},~\pi_1(Y)=H.$$
The splitting is {\it finite} if the complexes $W,X,Y$ are finite,
and {\it infinite} otherwise.

A connected CW complex $W$ with $\pi_1(W)=G=
\begin{cases} G_1*_HG_2\\
G_1*_H\{t\}
\end{cases}$ an injective generalized free product is a homotopy pushout
$$\begin{cases}
\xymatrix{\widetilde{W}/H \ar[d]_-{i_2}
\ar[r]^-{i_1} & \widetilde{W}/G_1 \ar[d]^-{j_1} \\
\widetilde{W}/G_2 \ar[r]^-{j_2} & W}\\[12ex]
\xymatrix{\widetilde{W}/H \times \{0,1\} \ar[d]
\ar[r]^-{i_1\cup i_2} & \widetilde{W}/G_1 \ar[d]^-{j_1} \\
\widetilde{W}/H \times [0,1] \ar[r] & W}
\end{cases}$$
with $\widetilde{W}$ the universal cover of $W$ and
$\begin{cases}
i_1,i_2,j_1,j_2\\
i_1,i_2,j_1
\end{cases}$ the covering projections.  (See Proposition $\begin{cases}
\ref{lemma}\\ \ref{lemma2} \end{cases}$ for proofs).
Thus $W$ has a canonical infinite injective Seifert--van Kampen splitting
$(X(\infty),Y(\infty))$ with
$$\begin{cases}
Y(\infty)=\widetilde{W}/H\times \{1/2\} \subset
X(\infty)=\widetilde{W}/G_1 \cup_{i_1} \widetilde{W}/H \times [0,1]\cup_{i_2}
\widetilde{W}/G_2\\[1ex]
Y(\infty)=\widetilde{W}/H\times \{1/2\} \subset
X(\infty)=\widetilde{W}/G_1 \cup_{i_1 \cup i_2}
\widetilde{W}/H \times [0,1]~~~.
\end{cases}$$
For finite $W$ with $\pi_1(W)$ a finitely presented injective generalized
free product it is easy to obtain finite injective Seifert--van Kampen
splittings by codimen\-sion--1 manifold transversality.  In fact, there
are two somewhat different ways of doing so:
\begin{itemize}
\item[(i)] Consider a regular neighbourhood $(M,\partial M)$ of $W
\subset S^n$ ($n$ large), apply codimension--1 manifold
transversality to a map
$$\begin{cases}
f\co M \to BG=BG_1\cup_{BH \times \{0\}}BH \times [0,1]
\cup_{BH \times\{1\}}BG_2 \\
f\co M \to BG=BG_1\cup_{BH \times
\{0,1\}}BH \times [0,1] \end{cases}$$
inducing the identification $\pi_1(M)=G$ to obtain a finite
Seifert--van Kampen splitting $(M^n,N^{n-1})$ with
$N=f^{-1}(BH\times \{1/2\}) \subset M$, and then make the splitting
injective by low-dimensional handle exchanges.
\item[(ii)] Replace the 2--skeleton $W^{(2)}$ by a homotopy equivalent
manifold with boundary $(M,\partial M)$, so $\pi_1(M)=\pi_1(W)$ is a
finitely presented injective generalized free product and $M$ has a
finite injective Seifert--van Kampen splitting by manifold transversality (as in
(i)).  Furthermore, $(W,M)$ is a finite CW pair and
$$W=M \cup \bigcup\limits_{n \geqslant 3}(W,M)^{(n)}$$
with the relative $n$--skeleton $(W,M)^{(n)}$ a union of $n$--cells $D^n$
attached along maps $S^{n-1} \to M \cup (W,M)^{(n-1)}$.  Set
$(W,M)^{(2)}=\emptyset$, and assume inductively that for some $n
\geqslant 3$ $M \cup (W,M)^{(n-1)}$ already has a finite Seifert--van Kampen
splitting $(X,Y)$.  For each $n$--cell $D^n \subset (W,M)^{(n)}$ use
manifold transversality to make the composite
$$S^{n-1} \to M \cup (W,M)^{(n-1)}\simeq X$$
transverse at $Y \subset X$, and extend this transversality to
make the composite
$$f\co D^n \to M \cup (W,M)^{(n)} \to BG$$
transverse at $BH \subset BG$.  The transversality gives $D^n$ a finite
CW structure in which $N^{n-1}=f^{-1}(BH) \subset D^n$ is a
subcomplex, and
$$(X',Y')=\bigg(X \cup \bigcup\limits_{D^n \subset (W,M)^{(n)}} D^n\,,\,
Y \cup\bigcup\limits_{D^n \subset (W,M)^{(n)}} N^{n-1}\bigg)$$
is an extension to $M \cup (W,M)^{(n)}$ of the finite Seifert--van Kampen
splitting.
\end{itemize}
However, the geometric nature of manifold transversality does not give
any insight into the CW structures of the splittings $(X,Y)$ of $W$
obtained as above, let alone into the algebraic analogue of
transversality for $\Z[G]$--module chain complexes.  Here, we obtain
Seifert--van Kampen splittings combinatorially, in the following
converse of the Seifert--van Kampen Theorem.

\begin{ctt}
Let $W$ be a finite connected CW complex with $\pi_1(W)=G=
\begin{cases} G_1*_HG_2\\
G_1*_H\{t\}
\end{cases}$ an injective generalized free product.

{\rm (i)}\qua The canonical infinite Seifert--van Kampen splitting
$(X(\infty),Y(\infty))$ of $W$ is a union of
finite Seifert--van Kampen splittings $(X,Y) \subset (X(\infty),Y(\infty))$
$$(X(\infty),Y(\infty))=\bigcup (X,Y).$$
In particular, there exist finite Seifert--van Kampen splittings
$(X,Y)$ of $W$.

{\rm (ii)}\qua If the injective generalized free product structure on
$\pi_1(W)$ is finitely presented then for any finite Seifert--van Kampen
splitting $(X,Y)$ of $W$ it is possible to attach finite numbers of 2--
and 3--cells to $X$ and $Y$ to obtain an injective finite Seifert--van
Kampen splitting $(X',Y')$ of $W$, such that $(X,Y) \subset (X',Y')$
with the inclusion $X \to X'$ a homotopy equivalence and the inclusion
$Y \to Y'$ a $\Z[H]$--coefficient homology equivalence.
\end{ctt}

The Theorem is proved in section \ref{comb}.  The main ingredient of the
proof is the construction of a finite Seifert--van Kampen splitting of
$W$ from a finite {\it domain} of the universal cover
$\widetilde{W}$, as given by finite subcomplexes
$\begin{cases}
W_1,W_2 \subseteq \widetilde{W}\\
W_1 \subseteq \widetilde{W}
\end{cases}$
such that
$$\begin{cases}
G_1W_1 \cup G_2 W_2 =\widetilde{W} \\
G_1W_1=\widetilde{W}~.
\end{cases}$$

Algebraic transversality makes much use of the induction and restriction
functors associated to a ring morphism $i\co A \to B$
$$\begin{array}{l}
i_!\co \{A\hbox{-modules}\} \to \{B\hbox{-modules}\};~
  M \mapsto i_!M=B\otimes_AM,\\[1ex]
i^!\co \{B\hbox{-modules}\} \to \{A\hbox{-modules}\};~
  N\mapsto i^!N=N.
\end{array}$$
These functors are adjoint, with
$${\rm Hom}_B(i_!M,N)={\rm Hom}_A(M,i^!N).$$
Let $G=\begin{cases} G_1*_HG_2 \\ G_1*_H\{t\} \end{cases}$
be a generalized free product.
By definition, a {\it Mayer--Vietoris splitting} (or
{\it presentation}) ${\mathcal E}$ of a $\Z[G]$--module chain complex $C$ is:
\begin{itemize}
\item[{\rm (A)}] an exact sequence of $\Z[G]$--module chain complexes
$${\mathcal E}\co 0 \to k_!D \xymatrix@C+15pt{\ar[r]^-{
\begin{pmatrix} 1 \otimes e_1 \\ 1 \otimes e_2 \end{pmatrix}}&}
(j_1)_!C_1 \oplus (j_2)_!C_2 \to C \to 0$$
with $C_1$ a $\Z[G_1]$--module chain complex, $C_2$ a $\Z[G_2]$--module
chain complex,  $D$ a $\Z[H]$--module chain complex,
$e_1\co (i_1)_!D \to C_1$ a $\Z[G_1]$--module chain map and
$e_2\co (i_2)_!D \to C_2$ a $\Z[G_2]$--module chain map,
\item[{\rm (B)}] an exact sequence of $\Z[G]$--module chain complexes
$${\mathcal E}\co 0 \to (j_1i_1)_!D
\xymatrix@C+45pt{\ar[r]^-{\displaystyle{1 \otimes e_1-t \otimes e_2}}&}
(j_1)_!C_1 \to C \to 0$$
with $C_1$ a $\Z[G_1]$--module chain complex, $D$ a $\Z[H]$--module chain
complex, and $e_1\co (i_1)_!D \to C_1$, $e_2\co (i_2)_!D \to C_1$
$\Z[G_1]$--module chain maps.
\end{itemize}
A Mayer--Vietoris splitting ${\mathcal E}$ is {\it finite} if every
chain complex in ${\mathcal E}$ is finite f.g.\  free, and {\it infinite}
otherwise. See section \ref{tree} for the construction of a (finite)
Mayer--Vietoris splitting of the cellular $\Z[\pi_1(X)]$--module
chain complex $C(\widetilde{X})$ of the universal cover $\widetilde{X}$
of a (finite) connected CW complex $X$ with a 2--sided connected
subcomplex $Y \subset X$ such that $\pi_1(Y) \to \pi_1(X)$ is injective.

For any injective generalized free product
$G=\begin{cases}
G_1*_HG_2 \\
G_1*_H\{t\} \end{cases}$
every free $\Z[G]$--module chain complex $C$
has a canonical infinite Mayer--Vietoris splitting
$$\begin{array}{ll}
{\rm (A)}&
{\mathcal E}(\infty)\co 0 \to k_!k^!C \to (j_1)_!j_1^!C \oplus
(j_2)_!j_2^!C \to C \to 0\\[1ex]
{\rm (B)}&
{\mathcal E}(\infty)\co 0 \to k_!k^!C \to (j_1)_!j_1^!C  \to C \to 0.
\end{array}$$
For finite $C$ we shall obtain finite Mayer--Vietoris splittings in
the following converse of the Mayer--Vietoris Theorem.

\begin{att}
Let $G=\begin{cases}
G_1*_HG_2 \\ G_1*_H\{t\} \end{cases}$ be an injective generalized
free product. For a finite f.g.\ free $\Z[G]$--module chain complex
$C$ the canonical infinite Mayer--Vietoris splitting ${\mathcal
E}(\infty)$ of $C$ is a union of finite Mayer--Vietoris splittings
${\mathcal E} \subset {\mathcal E}(\infty)$
$${\mathcal E}(\infty)=\bigcup {\mathcal E}.$$
In particular, there exist finite Mayer--Vietoris splittings
${\mathcal E}$ of $C$.
\end{att}

The existence of finite Mayer--Vietoris splittings was first proved by
Waldhausen \cite{wald1}, \cite{wald2}.  The proof of the Theorem in
section \ref{alg} is a simplification of the original argument,
using chain complex analogues of the CW domains.

Suppose now that $(X,Y)$ is the finite 2--sided CW pair defined by a
(compact) connected $n$--dimensional manifold $X^n$ together with a
connected codimension--1 submanifold $Y^{n-1} \subset X$ with trivial
normal bundle.  By definition, a homotopy equivalence $f\co M^n \to X$
from an $n$--dimensional manifold {\it splits} at $Y \subset X$ if $f$
is homotopic to a map (also denoted by $f$) which is transverse at $Y$,
such that the restriction $f\vert\co N^{n-1}=f^{-1}(Y) \to Y$ is also a
homotopy equivalence.  In general, homotopy equivalences do not split:
it is not possible to realize the Seifert--van Kampen splitting $X$ of
$M$ by a codimension--1 submanifold $N \subset M$.  For $(X,Y)$ with
injective $\pi_1(Y) \to \pi_1(X)$ there are algebraic $K$-- and
$L$--theory obstructions to splitting homotopy equivalences, involving
the Nil-groups of Waldhausen \cite{wald1}, \cite{wald2} and the
UNil-groups of Cappell \cite{cap1}, and for $n \geqslant 6$ these are
the complete obstructions to splitting.  As outlined in Ranicki
\cite[section 7.6]{ranicki1}, \cite[section 8]{ranicki2}, algebraic
transversality for chain complexes is an essential ingredient for
a systematic treatment of both the algebraic $K$-- and $L$--theory
obstructions.  The algebraic analogue of the combinatorial approach to
CW transversality worked out here will be used to provide such a treatment
in Ranicki \cite{ranicki5}.

Although the algebraic $K$-- and
$L$--theory of generalized free products will not actually be considered
here, it is worth noting that the early results of Higman
\cite{higman}, Bass, Heller and Swan \cite{bhs} and Stallings
\cite{stallings1} on the Whitehead groups of polynomial extensions and
free products were followed by the work of the dedicatee on the
Whitehead group of amalgamated free products (Casson \cite{casson})
prior to the general results of Waldhausen \cite{wald1}, \cite{wald2}
on the algebraic $K$--theory of generalized free products.

The algebraic $K$--theory spectrum $A(X)$ of a space (or simplicial set) $X$
was defined by Waldhausen \cite{wald3} to be the $K$--theory
 $$A(X)=K({\mathcal R}_f(X))$$
of the category ${\mathcal R}_f(X)$ of retractive spaces over $X$,
and also as
 $$A(X)=K(S \wedge G(X)_+)$$
with $S$ the sphere spectrum and $G(X)$ the loop group of $X$. See
H\"utte\-mann, Klein, Vogell, Waldhausen and Williams \cite{hkvww},
Schw\"anzl and Staffeldt \cite{schsta}, Sch\-w\"anzl, Staffeldt and
Waldhausen \cite{schstawal} for the current state of knowledge
concerning the Mayer--Vietoris-type decomposition of $A(X)$ for a
finite 2--sided CW pair $(X,Y)$. The $A$--theory splitting
theorems obtained there use the second form of the definition of
$A(X)$. The Combinatorial Transversality Theorem could perhaps be
used to obtain $A$--theory splitting theorems directly from the
first form of the definition, at least for injective $\pi_1(Y) \to
\pi_1(X)$.

I am grateful to Bob Edwards, Dirk Sch\"utz and the referee for useful
suggestions.

\section{The Seifert--van Kampen and Mayer--Vietoris Theorems}
\label{tree}

Following some standard material on covers and fundamental groups
we recall the well-known Bass--Serre theory relating
injective generalized free products and groups acting on trees.  The
Seifert--van Kampen theorem for the fundamental group $\pi_1(X)$ and the
Mayer--Vietoris theorem for the cellular $\Z[\pi_1(X)]$--module chain
complex $C(\widetilde{X})$ of the universal cover $\widetilde{X}$ of a
connected CW complex $X$ with a connected 2--sided subcomplex $Y
\subset X$ and injective $\pi_1(Y) \to \pi_1(X)$ are then deduced from the
construction of the universal cover $\widetilde{X}$ of $X$ by cutting
along $Y$, using the tree $T$ on which $\pi_1(X)$ acts.

\subsection{Covers}

Let $X$ be a connected CW complex with fundamental group $\pi_1(X)=G$
and universal covering projection $p\co \widetilde{X} \to X$,
with $G$ acting on the left of $\widetilde{X}$.
Let $C(\widetilde{X})$ be the cellular free (left) $\Z[G]$--module
chain complex. For any subgroup $H \subseteq G$ the covering
$Z=\widetilde{X}/H$ of $X$ has universal cover $\widetilde{Z}=\widetilde{X}$
with cellular $\Z[H]$--module chain complex
$$C(\widetilde{Z})=k^!C(\widetilde{X})$$
with $k\co \Z[H] \to \Z[G]$ the inclusion. For a connected subcomplex
$Y \subseteq X$ the inclusion $Y \to X$ induces an
injection $\pi_1(Y) \to \pi_1(X)=G$ if and only if the components
of $p^{-1}(Y) \subseteq \widetilde{X}$ are copies of the universal
cover $\widetilde{Y}$ of $Y$. Assuming this injectivity condition
we have
$$p^{-1}(Y)=\bigcup\limits_{g \in [G;H]}g\widetilde{Y}
\subset \widetilde{X}$$
with $H=\pi_1(Y) \subseteq G$ and $[G;H]$ the set of right $H$--cosets
$$g=xH \subseteq G~~(x \in G).$$
The cellular $\Z[G]$--module chain complex of $p^{-1}(Y)$ is induced
from the cellular $\Z[H]$--module chain complex of $\widetilde{Y}$
$$C(p^{-1}(Y))=k_!C(\widetilde{Y})=\Z[G]\otimes_{\Z[H]}C(\widetilde{Y})=
\bigoplus\limits_{g \in [G;H]}C(g\widetilde{Y})
\subseteq C(\widetilde{X}).$$
The inclusion $Y \to X$ of CW complexes induces an inclusion of
$\Z[H]$--module chain complexes
$$C(\widetilde{Y}) \to C(Z)=k^!C(\widetilde{X})$$
adjoint to the inclusion of $\Z[G]$--module chain complexes
$$C(p^{-1}(Y))=k_!C(\widetilde{Y}) \to C(\widetilde{X}).$$

\subsection{Amalgamated free products}

\begin{Theorem}[Serre \cite{serre}] \label{serre}
A group $G$ is {\rm (}isomorphic to{\rm )} an injective amalgamated
free product $G_1*_HG_2$ if and only if $G$ acts on a tree $T$ with
$$T/G=[0,1].$$
\end{Theorem}
{\bf Idea of proof}\qua Given an injective amalgamated free
product $G=G_1*_HG_2$ let $T$ be the tree defined by
$$T^{(0)}=[G;G_1] \cup [G;G_2],~T^{(1)}=[G;H].$$
The edge
$h \in [G;H]$ joins the unique vertices $g_1 \in [G;G_1]$, $g_2 \in [G;G_2]$ with
$$g_1 \cap g_2=h \subset G.$$
The group $G$ acts on $T$ by
$$G \times T \to T;~(g,x) \mapsto gx$$
with $T/G=[0,1]$.\\
Conversely, if a group $G$ acts on a tree $T$ with $T/G= [0,1]$ then
$G=G_1*_HG_2$ is an injective amalgamated free product with
$G_i\subseteq G$ the isotropy subgroup of $G_i\in T^{(0)}$ and $H
\subseteq G$ the isotropy subgroup of $H \in T^{(1)}$.
\endproof

If the amalgamated free product $G$ is nontrivial the tree $T$ is infinite.

\begin{Theorem} Let
$$X=X_1\cup_YX_2$$
be a connected CW complex which is a union of connected
subcomplexes such that the morphisms induced by the inclusions
$Y \to X_1$, $Y \to X_2$
$$i_1\co \pi_1(Y)=H \to \pi_1(X_1)=G_1,~~
i_2\co \pi_1(Y)=H \to \pi_1(X_2)=G_2$$
are injective, and let
$$G=G_1*_HG_2$$
with tree $T$.

{\rm (i)}\qua The universal cover $\widetilde{X}$ of $X$ is the union
of translates of the universal covers $\widetilde{X}_1,\widetilde{X}_2$
of $X_1,X_2$
$$\widetilde{X}=\bigcup\limits_{g_1 \in [G;G_1]}g_1\widetilde{X}_1
\cup_{\bigcup\limits_{h \in [G;H]}h\widetilde{Y}}
\bigcup\limits_{g_2 \in [G;G_2]}g_2\widetilde{X}_2.$$
with intersections translates of the universal cover $\widetilde{Y}$ of $Y$
$$g_1\widetilde{X}_1 \cap g_2\widetilde{X}_2=
\begin{cases}
h \widetilde{Y}&\hbox{if}~g_1 \cap g_2=h \in [G;H]\\
\emptyset&\hbox{otherwise}~.
\end{cases}$$
{\rm (ii)}\qua {\rm (Seifert--van Kampen)}\qua The fundamental group of
$X$ is the injective amalgamated free product
$$\pi_1(X)=G=G_1*_HG_2.$$
{\rm (iii)}\qua {\rm (Mayer--Vietoris)}\qua
The cellular $\Z[\pi_1(X)]$--module chain complex $C(\widetilde{X})$
has a Mayer--Vietoris splitting
$$\begin{array}{l}
0 \to k_!C(\widetilde{Y})
\xymatrix@C+20pt
{\ar[r]^-{\begin{pmatrix} 1\otimes e_1 \\ 1\otimes e_2 \end{pmatrix}}&}
(j_1)_!C(\widetilde{X}_1)\oplus (j_2)_!C(\widetilde{X}_2)\\[1ex]
\hskip175pt
\xymatrix@C+35pt{\ar[r]^-{\displaystyle{(f_1-f_2)}}&} C(\widetilde{X})\to 0
\end{array}$$
with $e_1\co Y \to X_1$, $e_2\co Y \to X_2$, $f_1\co X_1 \to X$, $f_2\co X_2 \to X$
the inclusions.
\end{Theorem}
\begin{proof} (i)\qua Consider first the special case
$G_1=G_2=H=\{1\}$. Every map $S^1 \to X=X_1\cup_YX_2$ is homotopic to one
which is transverse at $Y \subset X$ (also denoted $f$) with $f(0)=f(1)
\in Y$, so that $[0,1]$ can be decomposed as a union of closed intervals
$$[0,1]=\bigcup\limits^n_{i=0}[a_i,a_{i+1}]~~(0=a_0 < a_1< \dots < a_{n+1}=1)$$
with
$$f(a_i) \in Y,~f[a_i,a_{i+1}] \subseteq
\begin{cases}
X_1&{\rm if}~i~{\rm is~even}\\[1ex]
X_2&{\rm if}~i~{\rm is~odd}~.
\end{cases}$$
Choosing paths $g_i\co [0,1] \to Y$ joining $a_i$ to $a_{i+1}$ and using
$\pi_1(X_1)=\pi_1(X_2)=\{1\}$ on the loops $f\vert_{[a_i,a_{i+1}]} \cup
g_i\co S^1 \to
X_1$ (resp. $X_2$) for $i$ even (resp. odd) there is obtained a
contraction of $f\co S^1 \to X$, so that $\pi_1(X)=\{1\}$ and $X$ is its
own universal cover.

In the general case let
$$p_j\co \widetilde{X}_j \to X_j~(j=1,2),~q\co \widetilde{Y} \to Y$$
be the universal covering projections. Since $i_j\co H \to G_j$ is injective
$$(p_j)^{-1}(Y)=\bigcup\limits_{h_j \in [G_j;H]}h_j\widetilde{Y}.$$
The CW complex defined by
$$\widetilde{X}=\bigcup\limits_{g_1 \in [G;G_1]}g_1\widetilde{X}_1
\cup_{\bigcup\limits_{h \in [G;H]}h\widetilde{Y}}
\bigcup\limits_{g_2 \in [G;G_2]}g_2\widetilde{X}_2$$
is simply-connected by the special case, with a free $G$--action such that
$\widetilde{X}/G=X$,
so that $\widetilde{X}$ is the universal cover of $X$ and $\pi_1(X)=G$.

(ii)\qua The vertices of the tree $T$ correspond to the translates of $\widetilde{X}_1$,
$\widetilde{X}_2 \subset \widetilde{X}$, and the edges correspond to the
translates of $\widetilde{Y} \subset \widetilde{X}$. The free action of $G$
on $\widetilde{X}$ determines a (non-free) action of $G$ on $T$ with
$T/G=[0,1]$, and $\pi_1(X)=G=G_1*_HG_2$ by Theorem \ref{serre}.

(iii)\qua Immediate from the expression of $\widetilde{X}$ in (i) as a union of
copies of $\widetilde{X}_1$ and $\widetilde{X}_2$.
\end{proof}

Moreover, in the above situation there is defined a $G$--equivariant map
$\widetilde{f}\co \widetilde{X} \to T$ with quotient a map
$$f\co \widetilde{X}/G=X \to T/G=[0,1]$$
such that
$$X_1=f^{-1}([0,1/2]),~X_2=f^{-1}([1/2,1]),~Y=f^{-1}(1/2) \subset X.$$

\subsection{HNN extensions}

\begin{Theorem} \label{serre2}
A group $G$ is {\rm (}isomorphic to{\rm )} an injective HNN extension
$G_1*_H\{t\}$ if and only if $G$ acts on a tree $T$ with
$$T/G=S^1.$$
\end{Theorem}
{\bf Idea of proof}\qua Given an injective HNN extension $G=G_1*_H\{t\}$
let $T$ be the infinite tree defined by
$$T^{(0)}=[G;G_1],~T^{(1)}=[G;H],$$
identifying $H=i_1(H) \subseteq G$.
The edge $h \in [G;H]$ joins the unique vertices
$g_1,g_2 \in [G;G_1]$ with
$$g_1 \cap g_2t^{-1} = h \subset G.$$
The group $G$ acts on $T$ by
$$G \times T \to T;~(g,x) \mapsto gx$$
with $T/G=S^1$, $G_1 \subseteq G$
the isotropy subgroup of $G_1 \in T^{(0)}$ and $H \subseteq G$ the
isotropy subgroup of $H \in T^{(1)}$.

Conversely, if a group $G$ acts on a tree $T$ with $T/G=S^1$ then
$G=G_1*_H\{t\}$ is an injective HNN extension with $G_1 \subset G$
the isotropy group of $G_1 \in T^{(0)}$ and $H \subset G$ the isotropy
group of $H \in T^{(1)}$.
\endproof

\begin{Theorem} Let
$$X=X_1\cup_{Y \times \{0,1\}}Y \times [0,1]$$
be a connected CW complex which is a union of connected
subcomplexes such that the morphisms induced by the inclusions
$Y\times \{0\} \to X_1$, $Y \times \{1\} \to X_1$
$$i_1,i_2\co \pi_1(Y)=H \to \pi_1(X_1)=G_1$$
are injective, and let
$$G=G_1*_H\{t\}$$
with tree $T$.

{\rm (i)}\qua The universal cover $\widetilde{X}$ of $X$ is the union
of translates of the universal cover $\widetilde{X}_1$ of $X_1$
$$\widetilde{X}=\bigcup\limits_{g_1 \in [G:G_1]}g_1\widetilde{X}_1
\cup_{\bigcup\limits_{h \in [G_1;H]} (h\widetilde{Y}\cup ht\widetilde{Y})}
\bigcup\limits_{h \in [G_1;H]} h\widetilde{Y}\times [0,1]$$
with $\widetilde{Y}$ the universal cover $\widetilde{Y}$.

{\rm (ii)}\qua {\rm (Seifert--van Kampen)}\qua The fundamental group of $X$ is the
injective HNN extension
$$\pi_1(X)=G=G_1*_H\{t\}.$$
{\rm (iii)}\qua {\rm (Mayer--Vietoris)}\qua The cellular $\Z[\pi_1(X)]$--module
chain complex $C(\widetilde{X})$ has a Mayer--Vietoris splitting
$${\mathcal E}\co 0 \to k_!C(\widetilde{Y})
\xymatrix@C+35pt{\ar[r]^-{\displaystyle{1 \otimes e_1-t \otimes e_2}}&}
(j_1)_!C(\widetilde{X}_1) \xymatrix{\ar[r]^-{\displaystyle{f}_1}&}
C(\widetilde{X})\to 0$$
with $e_1,e_2\co Y \to X_1$, $f_1\co X_1 \to X$ the inclusions.
\end{Theorem}
\begin{proof} (i)\qua Consider first the special case $G_1=H=\{1\}$, so that
$G=\Z=\{t\}$. The projection $\widetilde{X} \to X$
is a simply-connected regular covering with group of covering translations
$\Z$, so that it is the universal covering of $X$ and $\pi_1(X)=\Z$.

In the general case let
$$p_1\co \widetilde{X}_1 \to X_1,~q\co \widetilde{Y} \to Y$$
be the universal covering projections. Since $i_j\co H \to G_1$ is injective
$$\begin{array}{l}
(p_1)^{-1}(Y\times \{0\})=
\bigcup\limits_{g_1 \in [G;H]}g_1\widetilde{Y},\\[1ex]
(p_1)^{-1}(Y\times \{1\})= \bigcup\limits_{g_2 \in
[G;tHt^{-1}]}g_2\widetilde{Y}.
\end{array}$$
The CW complex defined by
$\widetilde{X}=\bigcup\limits_{g_1 \in [G:G_1]}g_1\widetilde{X}_1$
is simply-connected and with a free $G$--action such that
$\widetilde{X}/G=X$,
so that $\widetilde{X}$ is the universal cover of $X$ and $\pi_1(X)=G$.

(ii)\qua The vertices of the tree $T$ correspond to the translates of $\widetilde{X}_1
\subset \widetilde{X}$, and the edges correspond to the
translates of $\widetilde{Y}\times [0,1] \subset \widetilde{X}$. The
free action of $G$
on $\widetilde{X}$ determines a (non-free) action of $G$ on $T$ with
$T/G=S^1$, and $\pi_1(X)=G=G_1*_H\{t\}$ by Theorem \ref{serre2}.

(iii)\qua It is immediate from the expression of $\widetilde{X}$ in (i) as a union of
copies of $\widetilde{X}_1$ that there is defined a short exact sequence
$$\begin{array}{l}
0 \to k_!C(\widetilde{Y})\oplus k_!C(\widetilde{Y})
\xymatrix@C+45pt
{\ar[r]^-{\begin{pmatrix} 1\otimes e_1 & t\otimes e_2 \\ 1 & 1 \end{pmatrix}}&}
(j_1)_!C(\widetilde{X}_1)\oplus k_!C(\widetilde{Y})\\[1ex]
\hskip188pt
\xymatrix@C+55pt{\ar[r]^-{\displaystyle{(f_1-f_1(1\otimes e_1))}}&} C(\widetilde{X})\to 0
\end{array}$$
which gives the Mayer--Vietoris splitting.
\end{proof}

Moreover, in the above situation there is defined a
$G$--equivariant map $\widetilde{f}\co \widetilde{X} \to T$ with quotient a map
$$f\co \widetilde{X}/G=X \to T/G=[0,1]/(0=1)=S^1$$
such that
$$X_1=f^{-1}[0,1/2],~Y\times [0,1]=f^{-1}[1/2,1] \subset X.$$

\section{Algebraic transversality} \label{alg}

We now investigate the algebraic transversality properties of
$\Z[G]$--module chain complexes, with $G$ an injective generalized free
product. The Algebraic Transversality Theorem stated in the Introduction
will now be proved, treating the cases of an amalgamated free
product and an HNN extension separately.

\subsection{Algebraic transversality for amalgamated free products}

Let
$$G=G_1*_HG_2$$
be an injective amalgamated free product.
As in the Introduction write the injections as
$$\begin{array}{l}
i_1\co H \to G_1,~i_2\co H \to G_2,\\[1ex]
j_1\co G_1 \to G,~j_2\co G_2 \to G,\\[1ex]
k=j_1i_1=j_2i_2\co H \to G.
\end{array}$$

\begin{Definition} {\rm
{\rm (i)}\qua A {\it domain} $(C_1,C_2)$
of a $\Z[G]$--module chain complex $C$ is a pair of subcomplexes
$(C_1 \subseteq j_1^!C,C_2 \subseteq j_2^!C)$ such that the
chain maps
$$\begin{array}{l}
e_1\co (i_1)_!(C_1 \cap C_2) \to C_1;~b_1 \otimes y_1 \mapsto b_1y_1,\\[1ex]
e_2\co (i_2)_!(C_1 \cap C_2) \to C_2;~b_2 \otimes y_2 \mapsto b_2y_2,\\[1ex]
f_1\co (j_1)_!C_1 \to C;~a_1 \otimes x_1 \mapsto a_1x_1,\\[1ex]
f_2\co (j_1)_!C_2 \to C;~a_2 \otimes x_2 \mapsto a_2x_2
\end{array}$$
fit into a Mayer--Vietoris splitting of $C$
$${\mathcal E}(C_1,C_2)\co 0 \to k_!(C_1 \cap C_2)
\xymatrix{\ar[r]^-{\displaystyle{e}}&} (j_1)_!C_1 \oplus (j_2)_!C_2
\xymatrix{\ar[r]^-{\displaystyle{f}}&} C \to 0$$
with $e=\begin{pmatrix} e_1 \\ e_2\end{pmatrix}$, $f=(f_1-f_2)$.\\
(ii)\qua A domain $(C_1,C_2)$ is {\it finite} if $C_i$
($i=1,2$) is a finite f.g.\ free $\Z[G_i]$--module chain complex,
$C_1 \cap C_2$ is a finite f.g.\ free $\Z[H]$--module chain complex,
and {\it infinite} otherwise.}
\end{Definition}

\begin{Proposition}
Every free $\Z[G]$--module chain complex $C$ has a canonical
infinite domain $(C_1,C_2)=(j_1^!C,j_2^!C)$ with
$$C_1 \cap C_2=k^!C,$$
so that $C$ has a canonical infinite Mayer--Vietoris splitting
$${\mathcal E}(\infty)={\mathcal E}(j_1^!C,j_2^!C)\co 0 \to k_!k^!C \to
(j_1)_!j_1^! C \oplus (j_2)_!j_2^!C \to C \to 0.$$
\end{Proposition}
\begin{proof} It is enough to consider the special case $C=\Z[G]$,
concentrated in degree 0. The pair
$$(C_1,C_2)=(j_1^!\Z[G],j_2^!\Z[G])=(\bigoplus\limits_{[G;G_1]}\Z[G_1],
\bigoplus\limits_{[G;G_2]}\Z[G_2])$$
is a canonical infinite domain for $C$,  with
$${\mathcal E}(\infty)={\mathcal E}(C_1,C_2)\co 0 \to k_!k^!\Z[G] \to
(j_1)_!j_1^! \Z[G] \oplus (j_2)_!j_2^!\Z[G] \to \Z[G] \to 0$$
the simplicial chain complex $\Delta(T \times G)=\Delta(T)\otimes_{\Z}\Z[G]$,
along with its augmentation to $H_0(T \times G)=\Z[G]$.
\end{proof}

\begin{Definition} {\rm
(i)\qua For a based f.g.\ free $\Z[G]$--module $B=\Z[G]^b$ and
a subtree $U \subseteq T$ define a domain for $B$ (regarded
as a chain complex concentrated in degree 0)
$$(B(U)_1,B(U)_2)=(\sum\limits_{U^{(0)}_1}\Z[G_1]^b,
\sum\limits_{U^{(0)}_2}\Z[G_2]^b)$$
with
$$\begin{array}{l}
U^{(0)}_1=U^{(0)} \cap [G;G_1],~U^{(0)}_2=U^{(0)} \cap [G;G_2],\\[1ex]
B(U)_1 \cap B(U)_2= \sum\limits_{U^{(1)}}\Z[H]^b.
\end{array}$$
The associated Mayer--Vietoris splitting of $B$ is the
subobject ${\mathcal E}(U) \subseteq {\mathcal E}(\infty)$ with
$${\mathcal E}(U)\co 0 \to k_!\sum\limits_{U^{(1)}} \Z[H]^b
\to (j_1)_!\sum\limits_{U^{(0)}_1}\Z[G_1]^b \oplus
(j_2)_!\sum\limits_{U^{(0)}_2}\Z[G_2]^b \to B \to 0$$
the simplicial chain complex $\Delta(U \times G)^b=\Delta(U)\otimes_{\Z}B$,
along with its augmentation to $H_0(U \times G)^b=B$.
If $U \subset T$ is finite then $(B(U)_1,B(U)_2)$ is a finite domain.

(ii)\qua Let $C$ be an $n$--dimensional based f.g.\ free $\Z[G]$--module chain
complex, with $C_r=\Z[G]^{c_r}$.
A sequence $U=\{U_n,U_{n-1},\dots,U_1,U_0\}$ of subtrees $U_r \subseteq T$
is {\it realized} by $C$ if the differentials $d_C\co C_r \to C_{r-1}$ are
such that
$$d(C_r(U_r)_i) \subseteq C_{r-1}(U_{r-1})_i~~(1 \leqslant r \leqslant n,~i=1,2),$$
so that there is defined a Mayer--Vietoris splitting of $C$
$${\mathcal E}(U)\co 0 \to k_!\sum\limits_{U^{(1)}} C(U)_1 \cap C(U)_2
\to (j_1)_!\sum\limits_{U^{(0)}_1}C(U)_1 \oplus
(j_2)_!\sum\limits_{U^{(0)}_2}C(U)_2 \to C \to 0$$
with $C(U)_i$ the free $\Z[G_i]$--module chain complex defined by
$$d_{C(U)}=d_C\vert\co (C(U)_i)_r=C_r(U_r)_i
\to (C(U)_i)_{r-1}=C_{r-1}(U_{r-1})_i.$$
The sequence $U$ is {\it finite} if each subtree $U_r \subseteq T$
is finite, in which case ${\mathcal E}(U)$ is finite.}
\end{Definition}

\begin{Proposition}  \label{p}
For a finite based f.g.\ free $\Z[G]$--module chain complex $C$
the canonical infinite domain is a union of finite domains
$$(j_1^!C,j_2^!C)=\bigcup\limits_U (C(U)_1,C(U)_2),$$
with $U$ running over all the finite sequences which are realized by $C$.
The canonical infinite Mayer--Vietoris splitting of $C$
is thus a union of finite Mayer--Vietoris splittings
$${\mathcal E}(\infty)=\bigcup\limits_U {\mathcal E}(U).$$
\end{Proposition}
\begin{proof} The proof is based on the following observations:
\begin{itemize}
\item[(a)]  for any subtrees $V \subseteq U \subseteq T$
$${\mathcal E}(V) \subseteq {\mathcal E}(U) \subseteq {\mathcal E}(T)={\mathcal E}(\infty)$$
\item[(b)] the infinite tree $T$ is a union
$$T=\bigcup U$$
of the finite subtrees $U \subset T$,
\item[(c)] for any finite subtrees $U,U' \subset T$ there exists a
finite subtree $U'' \subset T$ such that $U \subseteq U''$ and
$U' \subseteq U''$,
\item[(d)] for every $d \in \Z[G]$ the $\Z[G]$--module morphism
$$d\co \Z[G] \to \Z[G];~x \mapsto xd$$
is resolved by a morphism $d_*\co {\mathcal E}(T) \to {\mathcal E}(T)$
of infinite Mayer--Vietoris splittings,
and for any finite subtree $U \subset T$ there exists a finite subtree
$U' \subset T$ such that
$$d_*({\mathcal E}(U)) \subseteq {\mathcal E}(U')$$
and $d_*\vert\co {\mathcal E}(U) \to {\mathcal E}(U')$ is a resolution
of $d$ by a morphism of finite Mayer--Vietoris splittings (cf. Proposition
1.1 of Waldhausen \cite{wald1}).
\end{itemize}
Assume $C$ is $n$--dimensional, with $C_r=\Z[G]^{c_r}$.
Starting with any  finite subtree $U_n \subseteq T$ let
$$U=\{U_n,U_{n-1},\dots,U_1,U_0\}$$
be a sequence of  finite subtrees $U_r \subset T$
such that the f.g.\ free submodules
$$\begin{array}{l}
C_r(U)_1=\sum\limits_{U^{(0)}_{r,1}}\Z[G_1]^{c_r} \subset j_1^!C_r=
\sum\limits_{T^{(0)}_1}\Z[G_1]^{c_r},\\
C_r(U)_2=\sum\limits_{U^{(0)}_{r,2}}\Z[G_2]^{c_r} \subset j_2^!C_r=
\sum\limits_{T^{(0)}_2}\Z[G_2]^{c_r},\\
D(U)_r=\sum\limits_{U^{(1)}_r}\Z[H]^{c_r} \subset k^!C_r=
\sum\limits_{T^{(1)}}\Z[H]^{c_r}
\end{array}$$
define subcomplexes
$$C(U)_1 \subset j_1^!C,~C(U)_2 \subset j_2^!C,~D(U) \subset k^!C.$$
Then $(C(U)_1,\allowbreak C(U)_2)$ is a domain of $C$
with
$$C(U)_1\cap C(U)_2=D(U),$$
and $U$ is realized by $C$.
\end{proof}

\begin{Remark} {\rm
(i)\qua The existence of finite Mayer--Vietoris splittings
was first proved by Waldhausen
\cite{wald1},\cite{wald2}, using essentially the same method.
See Quinn \cite{quinn} for a proof using controlled algebra.
The construction of generalized free products by
noncommutative localization (cf. Ranicki \cite{ranicki4})
can be used to provide a different proof.

(ii)\qua The construction of the finite Mayer--Vietoris splittings
${\mathcal E}(U)$ in \ref{p} as subobjects of the universal
Mayer--Vietoris splitting ${\mathcal E}(T)={\mathcal E}(\infty)$
is taken from Remark 8.7 of Ranicki \cite{ranicki2}.}
\end{Remark}

This completes the proof of the Algebraic Transversality Theorem
for amalgamated free products.

\subsection{Algebraic transversality for HNN extensions}

The proof of algebraic transversality for HNN extensions proceeds
exactly as for amalgamated free products, so only the statements
will be given.

Let
$$G=G_1*_H\{t\}$$
be an injective HNN extension. As in the Introduction, write the injections as
$$i_1,i_2\co H \to G_1,~j\co G_1 \to G,~k=j_1i_1=j_1i_2\co G_1 \to G.$$

\begin{Definition} {\rm
{\rm (i)}\qua A {\it domain} $C_1$
of a $\Z[G]$--module chain complex $C$ is a subcomplex
$C_1 \subseteq j_1^!C$ such that the chain maps
$$\begin{array}{l}
e_1\co (i_1)_!(C_1 \cap tC_1) \to C_1;~b_1 \otimes y_1 \mapsto b_1y_1,\\[1ex]
e_2\co (i_2)_!(C_1 \cap tC_1) \to C_1;~b_2 \otimes y_2 \mapsto b_2t^{-1}y_2,\\[1ex]
f\co (j_1)_!C_1 \to C;~a \otimes x \mapsto ax
\end{array}$$
fit into a Mayer--Vietoris splitting of $C$
$${\mathcal E}(C_1)\co 0 \to k_!(C_1 \cap tC_1)
\xymatrix@C+30pt{\ar[r]^-{\displaystyle{1\otimes e_1 -t\otimes e_2}}&} (j_1)_!C_1
\xymatrix{\ar[r]^-{\displaystyle{f}}&} C \to 0.$$
(ii)\qua A domain $C_1$ is {\it finite} if $C_1$
is a finite f.g.\ free $\Z[G_1]$--module chain complex and
$C_1 \cap tC_1$ is a finite f.g.\ free $\Z[H]$--module chain complex.}
\end{Definition}

\begin{Proposition}
Every free $\Z[G]$--module chain complex $C$ has a canonical
infinite domain $C_1=j_1^!C$ with
$$C_1 \cap tC_1=k^!C_1,$$
so that $C$ has a canonical infinite Mayer--Vietoris splitting
$${\mathcal E}(\infty)={\mathcal E}(j_1^!C)\co 0 \to k_!k^!C \to
(j_1)_!j_1^! C  \to C \to 0.\eqno{\qed}$$
\end{Proposition}

\begin{Definition} {\rm For any subtree $U \subseteq T$
define a domain for $\Z[G]$
$$C(U)_1=\sum\limits_{U^{(0)}}\Z[G_1]$$
with
$$C(U)_1 \cap tC(U)_1= \sum\limits_{U^{(1)}}\Z[H].$$
The associated Mayer--Vietoris splitting of $\Z[G]$ is the
subobject ${\mathcal E}(U) \subseteq {\mathcal E}(\infty)$ with
$${\mathcal E}(U)\co 0 \to k_!\sum\limits_{U^{(1)}} \Z[H]
\to (j_1)_!\sum\limits_{U^{(0)}}\Z[G_1] \to \Z[G] \to 0.$$
If $U \subset T$ is finite then $C(U)_1$ is finite.}
\end{Definition}

\begin{Proposition}  \label{q}
For a finite f.g.\ free $\Z[G]$--module chain complex $C$
the canonical infinite domain is a union of finite domains
$$j_1^!C=\bigcup C_1.$$
The canonical infinite Mayer--Vietoris splitting of $C$
is thus a union of finite Mayer--Vietoris splittings
$${\mathcal E}(\infty)=\bigcup {\mathcal E}(C_1).\eqno{\qed}$$
\end{Proposition}

This completes the proof of the Algebraic Transversality Theorem
for HNN extensions.

\section{Combinatorial transversality} \label{comb}

We now investigate the algebraic transversality properties of
CW complexes $X$ with $\pi_1(X)=G$ an injective generalized free
product. The Combinatorial Transversality Theorem stated in the Introduction
will now be proved, treating the cases of an amalgamated free
product and an HNN extension separately.

\subsection{Mapping cylinders}

We review some basic mapping cylinder constructions.

The {\it mapping cylinder} of a map $e\co V \to W$ is the identification space
$${\mathcal M}(e)=(V \times [0,1] \cup W)/\{(x,1) \sim e(x)\,\vert\, x\in V\}$$
such that $V=V\times \{0\} \subset {\mathcal M}(e)$. As ever, the projection
$$p\co {\mathcal M}(e) \to W;~\begin{cases}
(x,s) \mapsto e(x)&\hbox{\rm for $x \in V$, $s \in [0,1]$}\\
y \mapsto y&\hbox{\rm for $y \in W$}
\end{cases}$$
is a homotopy equivalence.

If $e$ is a cellular map of CW complexes then ${\mathcal M}(e)$ is a
CW complex. The cellular chain complex $C({\mathcal M}(e))$
is the {\it algebraic mapping cylinder} of the induced chain map
$e\co C(V) \to C(W)$, with
$$\begin{array}{l}
d_{C({\mathcal M}(e))}=\begin{pmatrix}
d_{C(W)} & (-1)^re & 0 \\
0 & d_{C(V)} & 0 \\
0 & (-1)^{r-1} & d_{C(V)} \end{pmatrix}\co \\[4ex]
C({\mathcal M}(e))_r=C(W)_r \oplus C(V)_{r-1} \oplus C(V)_r \\[1ex]
\hskip50pt \to C({\mathcal M}(e))_{r-1}=C(W)_{r-1} \oplus C(V)_{r-2} \oplus
C(V)_{r-1}.
\end{array}$$
The chain equivalence $p\co C({\mathcal M}(e)) \to C(W)$ is given by
$$p=(1~0~e)\co C({\mathcal M}(e))_r=C(W)_r \oplus C(V)_{r-1} \oplus C(V)_r
\to C(W)_r.$$

The {\it double mapping cylinder} ${\mathcal M}(e_1,e_2)$ of maps
$e_1\co V \to W_1$, $e_2\co V \to W_2$ is the identification space
$$\begin{array}{ll}
{\mathcal M}(e_1,e_2)&={\mathcal M}(e_1) \cup_V {\mathcal M}(e_2)\\[1ex]
&= W_1 \cup_{e_1} V \times [0,1] \cup_{e_2} W_2\\[1ex]
&= (W_1 \cup V \times [0,1] \cup W_2)/\{(x,0) \sim e_1(x), (x,1) \sim
e_2(x) \,\vert\, x \in V\}.
\end{array}$$
Given a commutative square of spaces and maps
$$\xymatrix{ V  \ar[r]^-{e_1} \ar[d]_-{e_2} &
W_1 \ar[d]^-{f_1} \\
W_2 \ar[r]^-{f_2} & W}$$
define the map
$$f_1 \cup f_2\co {\mathcal M}(e_1,e_2) \to W;~
\begin{cases}
(x,s) \mapsto f_1e_1(x)=f_2e_2(x)~~(x \in V,~s\in [0,1])\\
y_i \mapsto f_i(y_i)~~(y_i \in W_i,~i=1,2)~.
\end{cases}$$
The square is a {\it homotopy pushout} if $f_1 \cup f_2\co {\mathcal M}(e_1,e_2)
\to W$ is a homotopy equivalence.

If $e_1\co V \to W_1$, $e_2\co V \to W_2$ are cellular maps of CW complexes
then ${\mathcal M}(e_1,e_2)$ is a
CW complex, such that cellular chain complex $C({\mathcal M}(e_1,e_2))$
is the algebraic mapping cone of the chain map
$$\begin{pmatrix} e_1 \\ e_2 \end{pmatrix}\co C(V) \to C(W_1)\oplus C(W_2)$$
with
$$\begin{array}{l}
d_{C({\mathcal M}(e_1,e_2))}=\begin{pmatrix} d_{C(W_1)} & (-1)^re_1 & 0 \\
0 & d_{C(V)} & 0 \\
0 & (-1)^re_2 & d_{C(W_2)} \end{pmatrix}\co \\[4ex]
C({\mathcal M}(e_1,e_2))_r=C(W_1)_r \oplus C(V)_{r-1} \oplus C(W_2)_r \\[1ex]
\hskip50pt \to C({\mathcal M}(e_1,e_2))_{r-1}=C(W_1)_{r-1} \oplus C(V)_{r-2} \oplus
C(W_2)_{r-1}.
\end{array}$$

\subsection{Combinatorial transversality for amalgamated free pro\-ducts}

In this section $W$ is a connected CW complex with fundamental group
an injective amalgamated free product
$$\pi_1(W)=G=G_1*_HG_2$$
with tree $T$. Let $\widetilde{W}$ be the universal cover of $W$, and let
$$\xymatrix{ \widetilde{W}/H \ar[r]^-{i_1} \ar[d]_-{i_2} &
\widetilde{W}/G_1 \ar[d]^-{j_1} \\
\widetilde{W}/G_2 \ar[r]^-{j_2} & W}$$
be the commutative square of covering projections.

\begin{Definition} {\rm
{\rm (i)}\qua Suppose given subcomplexes $W_1,W_2 \subseteq \widetilde{W}$ such that
$$G_1W_1=W_1,~G_2W_2=W_2$$
so that
$$H(W_1 \cap W_2)=W_1 \cap W_2 \subseteq \widetilde{W}.$$
Define a commutative square of CW complexes and cellular maps
$$\xymatrix@C-10pt@R-10pt{
(W_1\cap W_2)/H \ar[rr]^-{e_1} \ar[dd]_-{e_2} & &
W_1/G_1 \ar[dd]^-{f_1} \\
& \Phi & \\
W_2/G_2 \ar[rr]^-{f_2} && W}$$
with
$$\begin{array}{l}
(W_1 \cap W_2)/H \subseteq \widetilde{W}/H,~
W_1/G_1 \subseteq \widetilde{W}/G_1,~
W_2/G_2 \subseteq \widetilde{W}/G_2,\\[1ex]
e_1=i_1\vert\co (W_1 \cap W_2)/H \to W_1/G_1,~e_2=i_2\vert\co (W_1 \cap W_2)/H \to W_2/G_2,\\[1ex]
f_1=j_1\vert\co W_1/G_1 \to W,~f_2=j_2\vert\co W_2/G_2\to W.
\end{array}$$
{\rm (ii)}\qua A {\it domain} $(W_1,W_2)$ for the
universal cover $\widetilde{W}$ of $W$ consists of connected subcomplexes
$W_1,W_2 \subseteq \widetilde{W}$ such that $W_1 \cap W_2$ is connected,
and such that for each cell $D \subseteq \widetilde{W}$
the subgraph $U(D) \subseteq T$ defined by
$$\begin{array}{l}
U(D)^{(0)}=\{g_1 \in [G;G_1]\,\vert\,
g_1 D \subseteq  W_1\} \cup
\{g_2 \in [G;G_2]\,\vert\,g_2 D \subseteq  W_2\}\\[1ex]
U(D)^{(1)}=\{h \in [G;H]\,\vert\,
h D \subseteq  W_1 \cap W_2\}
\end{array}$$
is a  tree.

{\rm (iii)}\qua A domain $(W_1,W_2)$ for $\widetilde{W}$ is
{\it fundamental} if the subtrees $U(D) \subseteq T$ are
either single vertices or single edges, so that
$$\begin{array}{l}
g_1W_1 \cap g_2W_2=\begin{cases} h(W_1 \cap W_2)&
\hbox{if $g_1 \cap g_2=h \in [G;H]$}\\
\emptyset&\hbox{if $g_1 \cap g_2=\emptyset$},
\end{cases}\\[3ex]
W=(W_1/G_1)\cup_{(W_1 \cap W_2)/H}(W_2/G_2).
\end{array}$$}
\end{Definition}

\begin{Proposition}
For a domain $(W_1,W_2)$ of $\widetilde{W}$ the pair of
cellular chain complexes $(C(W_1),C(W_2))$ is a domain
of the cellular chain complex $C(\widetilde{W})$.
\end{Proposition}
\begin{proof} The union of $GW_1,GW_2 \subseteq \widetilde{W}$ is
$$GW_1 \cup GW_2=\widetilde{W}$$
since for any cell $D \subseteq \widetilde{W}$ there
either exists $g_1 \in [G;G_1]$  such that
$g_1D \subseteq W_1$
or $g_2 \in [G;G_2]$ such that $g_2D \subseteq W_2$.
The intersection of $GW_1,GW_2 \subseteq \widetilde{W}$ is
$$GW_1 \cap GW_2=G(W_1 \cap W_2) \subseteq \widetilde{W}.$$
The Mayer--Vietoris exact sequence of cellular $\Z[G]$--module chain complexes
$$0 \to C(GW_1 \cap GW_2) \to C(GW_1) \oplus C(GW_2) \to C(\widetilde{W}) \to 0$$
is the Mayer--Vietoris splitting of $C(\widetilde{W})$ associated to $(C(W_1),C(W_2))$
$$0 \to k_!C(W_1 \cap W_2) \to (j_1)_!C(W_1) \oplus
(j_2)_!C(W_2) \to C(\widetilde{W}) \to 0$$
with $C(W_1 \cap W_2)=C(W_1)\cap C(W_2)$.
\end{proof}

\begin{Example} {\rm $W$ has a canonical infinite
 domain $(W_1,W_2)=(\widetilde{W},\widetilde{W})$
with $(W_1\cap W_2)/H=\widetilde{W}/H$,
and $U(D)=T$ for each cell $D \subseteq \widetilde{W}$.}
\end{Example}

\begin{Example} {\rm
(i)\qua Suppose that $W=X_1\cup_YX_2$, with $X_1,X_2,Y \subseteq W$
connected subcomplexes such that the isomorphism
$$\pi_1(W)=\pi_1(X_1)*_{\pi_1(Y)}\pi_1(X_2)
\xymatrix{\ar[r]^-{\displaystyle{\cong}}&} G=G_1*_HG_2$$
preserves the amalgamated free structures. Thus
$(W,Y)$ is a Seifert--van Kampen splitting of $W$,
and the morphisms
$$\pi_1(X_1) \to G_1,~\pi_1(X_2) \to G_2,~\pi_1(Y) \to H$$
are surjective.
(If $\pi_1(Y) \to \pi_1(X_1)$ and $\pi_1(Y) \to \pi_1(X_2)$ are
injective these morphisms are isomorphisms, and the splitting is
injective). The universal cover of $W$ is
$$\widetilde{W}=\bigcup\limits_{g_1 \in [G;G_1]}g_1\widetilde{X}_1
\cup_{\bigcup\limits_{h \in [G;H]}h\widetilde{Y}}
\bigcup\limits_{g_2 \in [G;G_2]}g_2\widetilde{X}_2$$
with $\widetilde{X}_i$ the regular cover of $X_i$ corresponding to
${\rm ker}(\pi_1(X_i) \to G_i)$ ($i=1,2$) and
$\widetilde{Y}$ the regular cover of $Y$ corresponding to
${\rm ker}(\pi_1(Y) \to H)$ (which are the universal covers
of $X_1,X_2,Y$ in the case
$\pi_1(X_1)=G_1$, $\pi_1(X_2)=G_2$, $\pi_1(Y)=H$).
The pair
$$(W_1,W_2)=(\widetilde{X}_1,\widetilde{X}_2)$$
is a fundamental domain of $\widetilde{W}$ such that
$$\begin{array}{l}
(W_1 \cap W_2)/H=Y,\\[1ex]
g_1W_1 \cap g_2W_2=
(g_1 \cap g_2)\widetilde{Y} \subseteq \widetilde{W}~~(g_1 \in [G;G_1], g_2 \in [G;G_2]).
\end{array}$$
For any cell $D \subseteq \widetilde{W}$
$$U(D)=\begin{cases}
\{g_1\}&\hbox{if $g_1D \subseteq
\widetilde{X}_1-\bigcup\limits_{h_1 \in [G_1;H]}h_1\widetilde{Y}$
for some $g_1 \in [G;G_1]$}\\
\{g_2\}&\hbox{if $g_2D \subseteq
\widetilde{X}_2-\bigcup\limits_{h_2 \in [G_2;H]}h_2\widetilde{Y}$
for some $g_2 \in [G;G_1]$}\\
\{g_1,g_2,h\}&\hbox{if $hD \subseteq \widetilde{Y}$ for
some $h=g_1\cap g_2 \in [G;H]$.}
\end{cases}$$
(ii)\qua If $(W_1,W_2)$ is a fundamental domain for any connected
CW complex $W$ with $\pi_1(W)=G=G_1*_HG_2$ then $W=X_1\cup_YX_2$
as in (i), with
$$X_1=W_1/G_1,~X_2=W_2/G_2,~Y=(W_1 \cap W_2)/H.$$}
\end{Example}

\begin{Definition} {\rm Suppose that $W$ is $n$--dimensional.
Lift each cell $D^r \subseteq W$ to a cell
$\widetilde{D}^r \subseteq \widetilde{W}$.
A sequence $U=\{U_n,U_{n-1},\dots,U_1,U_0\}$ of subtrees
$U_r \subseteq T$ is {\it realized} by $W$ if the subspaces
$$W(U)_1=\bigcup\limits^n_{r=0}\bigcup\limits_{D^r \subset W}
\bigcup\limits_{g_1 \in U^{(0)}_{r,1}}g_1\widetilde{D}^r,~
W(U)_2=\bigcup\limits^n_{r=0}\bigcup
\limits_{D^r \subset W}\bigcup\limits_{g_2 \in U^{(0)}_{r,2}}g_2\widetilde{D}^r
\subseteq \widetilde{W}$$
are connected subcomplexes, in which case $(W(U)_1,W(U)_2)$
is a domain for $\widetilde{W}$ with
$$W(U)_1\cap W(U)_2=
\bigcup\limits^n_{r=0}\bigcup\limits_{D^r \subset W}
\bigcup\limits_{h \in U^{(1)}_r}h\widetilde{D}^r \subseteq \widetilde{W}$$
a connected subcomplex. Thus $U$ is realized by $C(\widetilde{W})$
and
$$(C(W(U)_1),C(W(U)_2))=(C(\widetilde{W})(U)_1,C(\widetilde{W})(U)_2)
\subseteq (C(\widetilde{W}),C(\widetilde{W}))$$
is the domain for $C(\widetilde{W})$ given by
$(C_r(\widetilde{W})_1(U_r),C_r(\widetilde{W})(U)_2)$
in degree $r$.}
\end{Definition}

If a sequence $U=\{U_n,U_{n-1},\dots,U_1,U_0\}$ realized by $W$
is finite (ie if each $U_r\subseteq T$ is a finite subtree)
then $(W(U)_1,W(U)_2)$ is a finite domain for $\widetilde{W}$.

\begin{Proposition} \label{lemma}
{\rm (i)}\qua For any domain $(W_1,W_2)$ there is defined a
homotopy pushout
$$\xymatrix@C-10pt@R-10pt{
(W_1\cap W_2)/H \ar[rr]^-{e_1} \ar[dd]_-{e_2} & &
W_1/G_1 \ar[dd]^-{f_1} \\
& \Phi & \\
W_2/G_2 \ar[rr]^-{f_2} && W}$$
with $e_1=i_1\vert$, $e_2=i_2\vert$, $f_1=j_1\vert$, $f_2=j_2\vert$.
The connected 2--sided CW pair
$$(X,Y)=({\mathcal M}(e_1,e_2),(W_1 \cap W_2)/H \times \{1/2\})$$
is a Seifert--van Kampen splitting of $W$, with a homotopy equivalence
$$f=f_1\cup f_2\co X={\mathcal M}(e_1,e_2)
\xymatrix{\ar[r]^-{\displaystyle{\simeq}}&} W.$$
{\rm (ii)}\qua The commutative square of covering projections
$$\xymatrix{ \widetilde{W}/H \ar[r]^-{i_1} \ar[d]_-{i_2} &
\widetilde{W}/G_1 \ar[d]^-{j_1} \\
\widetilde{W}/G_2 \ar[r]^-{j_2} & W}$$
is a homotopy pushout. The connected 2--sided CW pair
$$(X(\infty),Y(\infty))=({\mathcal M}(i_1,i_2),\widetilde{W}/H \times \{1/2\})$$
is a canonical injective infinite Seifert--van Kampen splitting of $W$,
with a homotopy equivalence $j=j_1\cup j_2\co X(\infty) \to W$ such that
$$\pi_1(Y(\infty))=H \subseteq \pi_1(X(\infty))=G_1*_HG_2.$$
{\rm (iii)}\qua For any (finite) sequence $U=\{U_n,U_{n-1},\dots,U_0\}$
of subtrees of $T$ realized by $W$
there is defined a homotopy pushout
$$\xymatrix{ Y(U) \ar[r]^-{e_1} \ar[d]_-{e_2} &
X(U)_1 \ar[d]^-{f_1} \\
X(U)_2 \ar[r]^-{f_2} & W}$$
with
$$\begin{array}{l}
X(U)_1=W(U)_1/G_1,~X(U)_2=W(U)_2/G_2,\\[1ex]
Y(U)=(W(U)_1 \cap W(U)_2)/H,\\[1ex]
e_1=i_1\vert,~e_2=i_2\vert,~f_1=j_1\vert,~f_2=j_2\vert.
\end{array}$$
Thus
$$(X(U),Y(U))=({\mathcal M}(e_1,e_2),Y(U) \times \{1/2\})$$
is a (finite) Seifert--van Kampen splitting of $W$.

{\rm (iv)}\qua The canonical infinite domain of a finite CW complex $W$
with $\pi_1(W)=G_1*_HG_2$ is a union of finite domains
$$(\widetilde{W},\widetilde{W})=\bigcup\limits_U(W(U)_1,W(U)_2)$$
with $U$ running over all the finite sequences realized by  $W$.
The canonical infinite Seifert--van Kampen splitting of $W$ is thus a
union of finite Seifert--van Kampen splittings
$$(X(\infty),Y(\infty))=\bigcup\limits_U(X(U),Y(U)).$$
\end{Proposition}
\proof
(i)\qua Given a cell $D \subseteq W$ let
$\widetilde{D} \subseteq \widetilde{W}$ be a lift.
The inverse image of the interior ${\rm int}(D) \subseteq W$
$$f^{-1}({\rm int}(D))=U(\widetilde{D})\times {\rm int}(D) \subseteq
{\mathcal M}(i_1,i_2)=T\times_G \widetilde{W}$$
is contractible. In particular, point inverses are contractible,
so that $f\co X \to W$ is a homotopy equivalence. (Here is a more direct
proof that $f\co X \to W$ is a $\Z[G]$--coefficient homology equivalence.
The Mayer--Vietoris Theorem applied to the union
$\widetilde{W}=GW_1\cup GW_2$
expresses $C(\widetilde{W})$ as the cokernel of the $\Z[G]$--module chain map
$$e=\begin{pmatrix} 1 \otimes e_1 \\ 1 \otimes e_2 \end{pmatrix}\co
\Z[G]\otimes_{\Z[H]}C(W_1 \cap W_2) \to
\Z[G]\otimes_{\Z[G_1]}C(W_1) \oplus \Z[G]\otimes_{\Z[G_2]}C(W_2)$$
with a Mayer--Vietoris splitting
$$\begin{array}{l}
0 \to \Z[G]\otimes_{\Z[H]}C(W_1 \cap W_2)
\xymatrix{\ar[r]^-{\displaystyle{e}}&} \Z[G]\otimes_{\Z[G_1]}C(W_1) \oplus
Z[G]\otimes_{\Z[G_2]}C(W_2)\\[1ex]
\hskip250pt \xymatrix{\ar[r]&} C(\widetilde{W}) \to 0.
\end{array}$$
The decomposition $X={\mathcal M}(e_1,e_2)=X_1\cup_YX_2$ with
$$X_i={\mathcal M}(e_i)~(i=1,2),~
Y=X_1\cap X_2=(W_1 \cap W_2)/H\times \{1/2\}$$
lifts to a decomposition of the universal cover as
$$\widetilde{X}=\bigcup\limits_{g_1 \in [G;G_1]}g_1\widetilde{X}_1
\cup_{\bigcup\limits_{h \in [G;H]}h\widetilde{Y}}
\bigcup\limits_{g_2 \in [G;G_2]}g_2\widetilde{X}_2.$$
The Mayer--Vietoris splitting
$$0 \to \Z[G]\otimes_{\Z[H]}C(\widetilde{Y})
\xymatrix{\ar[r]&} \Z[G]\otimes_{\Z[G_1]}C(\widetilde{X}_1) \oplus
\Z[G]\otimes_{\Z[G_2]}C(\widetilde{X}_2) \!\to\! C(\widetilde{X}) \to 0,$$
expresses $C(\widetilde{X})$ as the algebraic mapping cone of the chain map $e$
$$C(\widetilde{X}) = {\mathcal C}(e\co \Z[G]\otimes_{\Z[H]}C(W_1 \cap W_2)
\to \Z[G]\otimes_{\Z[G_1]}C(W_1) \oplus \Z[G]\otimes_{\Z[G_2]}C(W_2)).$$
Since $e$ is injective the $\Z[G]$--module chain map
$$\widetilde{f}={\rm projection}\co C(\widetilde{X})={\mathcal C}(e)  \to
C(\widetilde{W})={\rm coker}(e)$$
induces isomorphisms in homology.)

(ii)\qua Apply (i) to $(W_1,W_2)=(\widetilde{W},\widetilde{W})$.

(iii)\qua Apply (i) to the domain $(W(U)_1,W(U)_2)$.

(iv)\qua Assume that $W$ is $n$--dimensional.  Proceed as for the chain
complex case in the proof of Proposition \ref{p} for the existence of a
domain for $C(\widetilde{W})$, but use only the sequences
$U=\{U_n,U_{n-1},\dots,U_0\}$ of finite subtrees $U_r \subset T$
realized by $W$.  An arbitrary finite subtree $U_n \subset T$ extends
to a finite sequence $U$ realized by $W$ since for $r \geqslant 2$ each
$r$--cell $\widetilde{D}^r \subset \widetilde{W}$ is attached to an
$(r-1)$--dimensional finite connected subcomplex, and every 1--cell
$\widetilde{D}^1 \subset \widetilde{W}$ is contained in a 1--dimensional
finite connected subcomplex.  Thus finite sequences $U$ realized by $W$
exist, and can be chosen to contain arbitrary finite collections of
cells of $\widetilde{W}$, with
$$(\widetilde{W},\widetilde{W})=\bigcup\limits_U(W(U)_1,W(U)_2).\eqno{\qed}$$

This completes the proof of part (i) of the Combinatorial
Transversality Theorem, the existence of finite Seifert--van Kampen
splittings.  Part (ii) deals with existence of finite injective
Seifert--van Kampen splittings: the adjustment of fundamental groups
needed to replace $(X(U),Y(U))$ by a homology-equivalent finite
injective Seifert--van Kampen splitting will use the following
rudimentary version of the Quillen plus construction.

\begin{Lemma} \label{l}
Let $K$ be a connected CW complex with a finitely
generated fundamental group $\pi_1(K)$.
For any surjection $\phi\co \pi_1(K) \to \Pi$
onto a finitely presented group $\Pi$ it is possible to attach a finite
number $n$ of 2-- and 3--cells to $K$ to obtain a connected CW complex
$$K'=K \cup \bigcup\limits_n D^2 \cup \bigcup\limits_n D^3$$
such that the inclusion $K \to K'$ is a $\Z[\Pi]$--coefficient
homology equivalence inducing $\phi\co \pi_1(K) \to \pi_1(K')=\Pi$.
\end{Lemma}
\begin{proof} The kernel of $\phi\co \pi_1(K) \to \Pi$ is the normal closure of a
finitely generated subgroup $N \subseteq \pi_1(K)$ by Lemma I.4 of Cappell \cite{cap2}.
(Here is the proof. Choose finite generating sets
$$g=\{g_1,g_2,\dots,g_r\} \subseteq \pi_1(K),~
h=\{h_1,h_2,\dots,h_s\} \subseteq \Pi$$
and let $w_k(h_1,h_2,\dots,h_s)$
($1 \leqslant  k \leqslant t$) be words in $h$ which are relations for $\Pi$.
As $\phi$ is surjective, can choose $h'_j \in \pi_1(K)$ with
$\phi(h'_j)=h_j$ ($1 \leqslant j \leqslant s$).
As $h$ generates $\Pi$ $\phi(g_i)=v_i(h_1,h_2,\dots,h_s)$ $(1 \leqslant
i \leqslant r)$ for some words $v_i$ in $h$. The kernel of $\phi$ is the
normal closure $N=\langle N'\rangle\triangleleft \pi_1(K)$ of the
subgroup $N'\subseteq \pi_1(K)$ generated by the finite set
$\{v_i(h'_1,\dots,h'_s)g_i^{-1},w_k(h'_1,\dots,h'_s)\}$.)
Let $x=\{x_1,x_2,\dots,x_n\}\subseteq \pi_1(K)$ be a finite set of generators of $N$,
and set
$$L=K \cup_x \bigcup\limits^n_{i=1} D^2.$$
The inclusion $K \to L$ induces
$$\phi\co \pi_1(K) \to \pi_1(L)=\pi_1(K)/\langle x_1,x_2,\dots,x_n \rangle=
\pi_1(K)/\langle N \rangle=\Pi.$$
Let $\widetilde{L}$ be the universal cover of $L$, and let $\widetilde{K}$
be the pullback cover of $K$. Now
$$\pi_1(\widetilde{K})={\rm ker}(\phi)=\langle x_1,x_2,\dots,x_n \rangle=
\langle N \rangle$$
so that the attaching maps $x_i\co S^1 \to K$ of the
2--cells in $L-K$ lift to null-homotopic maps $\widetilde{x}_i\co
S^1 \to \widetilde{K}$. The cellular chain complexes of $\widetilde{K}$
and $\widetilde{L}$ are related by
$$C(\widetilde{L})=C(\widetilde{K}) \oplus \bigoplus\limits_n (\Z[\Pi],2)$$
where $(\Z[\Pi],2)$ is just $\Z[\Pi]$ concentrated in degree 2. Define
$$x^*=\{x^*_1,x^*_2,\dots,x^*_n\}\subseteq \pi_2(L)$$
by
$$x^*_i = (0,(0,\dots,0,1,0,\dots,0)) \in \pi_2(L) = H_2(\widetilde{L})
 = H_2(\widetilde{K}) \oplus \Z[\Pi]^n~(1 \leqslant i \leqslant n),$$
and set
$$K'=L \cup_{x^*} \bigcup\limits^n_{i=1} D^3.$$
The inclusion $K \to K'$ induces $\phi\co\pi_1(K) \to \pi_1(K')=\pi_1(L)=\Pi$,
and the relative cellular $\Z[\Pi]$--module chain complex is
$$C(\widetilde{K'},\widetilde{K})\co \dots \to 0 \to \Z[\Pi]^n
\xymatrix{\ar[r]^-{1}&} \Z[\Pi]^n \to 0 \to \cdots$$
concentrated in degrees 2,3. In particular, $K \to K'$
is a $\Z[\Pi]$--coefficient homology equivalence.
\end{proof}

\begin{Proposition}
Let $(X,Y)$ be a finite connected 2--sided CW pair
with $X=X_1\cup_YX_2$ for connected $X_1,X_2,Y$,
together with an isomorphism
$$\pi_1(X)=\pi_1(X_1)*_{\pi_1(Y)}\pi_1(X_2)
\xymatrix{\ar[r]^-{\displaystyle{\cong}}&} G=G_1*_HG_2$$
preserving amalgamated free product structures, with the structure on $G$
injective. It is possible to attach a finite number of
2-- and 3--cells to $(X,Y)$ to obtain a finite injective
Seifert--van Kampen splitting $(X',Y')$ with
$X'=X'_1\cup_{Y'}X'_2$ such that
\begin{itemize}
\item[\rm (i)] $\pi_1(X')=G$, $\pi_1(X'_i)=G_i$ $(i=1,2)$, $\pi_1(Y')=H$,
\item[\rm (ii)] the inclusion $X \to X'$ is a homotopy equivalence,
\item[\rm (iii)] the inclusion $X_i \to X_i'$  $(i=1,2)$
is a $\Z[G_i]$--coefficient homology equivalence,
\item[\rm (iv)] the inclusion $Y \to Y'$ is a $\Z[H]$--coefficient
homology equivalence.
\end{itemize}
\end{Proposition}
\begin{proof} Apply the construction of Lemma \ref{l} to the
surjections $\pi_1(X_1) \to G_1$, $\pi_1(X_2) \to G_2$,
$\pi_1(Y) \to H$, to obtain
$$\begin{array}{l}
X'_i=(X_i \cup_{x_i}\bigcup\limits_{m_i}D^2)
\cup_{x^*_i}\bigcup\limits_{m_i}D^3~(i=1,2),\\[2ex]
Y'=(Y \cup_y\bigcup\limits_nD^2)\cup_{y^*}\bigcup\limits_nD^3
\end{array}
$$
for any $y=\{y_1,y_2,\dots,y_n\}\subseteq\pi_1(Y)$
such that ${\rm ker}(\pi_1(Y) \to H)$ is the normal closure of the subgroup
of $\pi_1(Y)$ generated by $y$, and any
$$x_i=\{x_{i,1},x_{i,2},\dots,x_{i,m_i}\}\subseteq\pi_1(X_i)$$
such that ${\rm ker}(\pi_1(X_i) \to G_i)$ is the normal closure of the
subgroup of $\pi_1(X_i)$ generated by $x_i$ ($i=1,2$).
Choosing $x_1,x_2$ to contain the images of $y$, we obtain the required
2--sided CW pair $(X',Y')$ with $X'=X'_1\cup_{Y'}X'_2$.
\end{proof}

This completes the proof of the Combinatorial Transversality Theorem
for amalgamated free products.

\subsection{Combinatorial transversality for HNN extensions}

The proof of combinatorial transversality for HNN extensions proceeds
exactly as for amalgamated free products, so only the statements will
be given.

In this section $W$ is a connected CW complex with fundamental group
an injective HNN extension
$$\pi_1(W)=G=G_1*_H\{t\}$$
with tree $T$. Let $\widetilde{W}$ be the universal cover of $W$, and let
$$\widetilde{W}/H
\raise5pt\hbox{\xymatrix@R-25pt@C-5pt{\ar[r]^-{i_1}&\\ \ar[r]_{i_2}&}}
\xymatrix{ \widetilde{W}/G_1 \ar[r]^-{j_1} & W}$$
be the covering projections, and define a commutative square
$$\xymatrix{ \widetilde{W}/H \times \{0,1\}
\ar[r]^-{i_1\cup i_2} \ar[d]_-{i_3} &
\widetilde{W}/G_1 \ar[d]^-{j_1} \\
\widetilde{W}/H \times [0,1] \ar[r]^-{j_2} & W}$$
where
$$\begin{array}{l}
i_3={\rm inclusion}\co \widetilde{W}/H \times \{0,1\} \to
\widetilde{W}/H \times [0,1],\\[1ex]
j_2\co \widetilde{W}/H \times [0,1] \to  W;~(x,s) \mapsto j_1i_1(x)=j_1i_2(x).
\end{array}$$

\begin{Definition} {\rm
{\rm (i)}\qua Suppose given a subcomplex $W_1 \subseteq \widetilde{W}$ with
$$G_1W_1=W_1$$
so that
$$H(W_1 \cap tW_1)=W_1 \cap tW_1 \subseteq \widetilde{W}.$$
Define a commutative square of CW complexes and cellular maps
$$\xymatrix@C-10pt@R-10pt{
(W_1\cap tW_1)/H \times \{0,1\}
\ar[rr]^-{e_1} \ar[dd]_-{e_2} & &
W_1/G_1 \ar[dd]^-{f_1} \\
& \Phi & \\
(W_1\cap tW_1)/H \times [0,1] \ar[rr]^-{f_2} && W}$$
with
$$\begin{array}{l}
(W_1 \cap tW_1)/H \subseteq \widetilde{W}/H,~W_1/G_1 \subseteq \widetilde{W}/G_1,\\[1ex]
e_1=(i_1 \cup i_2)\vert\co (W_1 \cap tW_1)/H \times \{0,1\} \to W_1/G_1,\\[1ex]
e_2=i_3 \vert\co
(W_1 \cap tW_1)/H \times \{0,1\} \to (W_1 \cap tW_1)/H \times [0,1],\\[1ex]
f_1=j_1\vert\co W_1/G_1 \to W,~f_2=j_2\vert\co (W_1\cap tW_1)/H \times [0,1] \to W.
\end{array}$$
{\rm (ii)}\qua A {\it domain} $W_1$ for the
universal cover $\widetilde{W}$ of $W$ is a connected subcomplex
$W_1 \subseteq \widetilde{W}$ such that $W_1 \cap tW_1$ is connected,
and such that for each cell $D \subseteq \widetilde{W}$
the subgraph $U(D) \subseteq T$ defined by
$$\begin{array}{l}
U(D)^{(0)}=\{g_1 \in [G;G_1]\,\vert\,g_1 D \subseteq  W_1\} \\[1ex]
U(D)^{(1)}=\{h \in [G_1;H]\,\vert\, h D \subseteq  W_1 \cap tW_1\}
\end{array}$$
is a  tree.

{\rm (iii)}\qua A domain $W_1$ for $\widetilde{W}$ is
{\it fundamental} if the subtrees $U(D) \subseteq T$ are
either single vertices or single edges, so that
$$\begin{array}{l}
g_1W_1 \cap g_2W_1=\begin{cases} h(W_1 \cap tW_1)&
\hbox{if $g_1 \cap g_2t^{-1}=h \in [G_1;H]$}\\
g_1W_1&\hbox{if $g_1=g_2$}\\
\emptyset&\hbox{if $g_1 \neq g_2$ and $g_1 \cap g_2t^{-1}=\emptyset$},
\end{cases}\\[5ex]
W=(W_1/G_1)\cup_{(W_1 \cap tW_1)/H \times \{0,1\}}
(W_1 \cap tW_1)/H \times [0,1].
\end{array}$$}
\end{Definition}

\begin{Proposition}
For a domain $W_1$ of $\widetilde{W}$ the
cellular chain complex $C(W_1)$ is a domain
of the cellular chain complex $C(\widetilde{W})$.\hfill$\qed$
\end{Proposition}

\begin{Example} {\rm $W$ has a canonical infinite domain
$W_1=\widetilde{W}$ with
$$(W_1\cap tW_1)/H=\widetilde{W}/H$$
and $U(D)=T$ for each cell $D \subseteq \widetilde{W}$.}
\end{Example}

\begin{Example} {\rm
(i)\qua Suppose that $W=X_1\cup_{Y \times \{0,1\}}Y \times [0,1]$,
with $X_1,Y \subseteq W$ connected subcomplexes such that the isomorphism
$$\pi_1(W)=\pi_1(X_1)*_{\pi_1(Y)}\{t\} \xymatrix{\ar[r]^-{\displaystyle{\cong}}&} G=G_1*_H\{t\}$$
preserves the HNN extensions. The morphisms
$\pi_1(X_1) \to G_1$, $\pi_1(Y) \to H$
are surjective. (If $i_1,i_2\co \pi_1(Y) \to \pi_1(X_1)$ are
injective these morphisms are also injective, allowing identifications
$\pi_1(X_1)=G_1$, $\pi_1(Y)=H$).
The universal cover of $W$ is
$$\widetilde{W}=\bigcup\limits_{g_1 \in [G:G_1]}g_1\widetilde{X}_1
\cup_{\bigcup\limits_{h \in [G_1;H]} (h\widetilde{Y}\cup ht\widetilde{Y})}
\bigcup\limits_{h \in [G_1;H]} h\widetilde{Y}\times [0,1]$$
with $\widetilde{X}_1$ the regular cover of $X_1$ corresponding to
${\rm ker}(\pi_1(X_1) \to G_1)$  and
$\widetilde{Y}$ the regular cover of $Y$ corresponding to
${\rm ker}(\pi_1(Y) \to H)$ (which are the universal covers
of $X_1,Y$ in the case $\pi_1(X_1)=G_1$, $\pi_1(Y)=H$).
Then $W_1=\widetilde{X}_1$ is a fundamental domain of
$\widetilde{W}$ such that
$$\begin{array}{l}
(W_1 \cap tW_1)/H=Y,~W_1 \cap tW_1=\widetilde{Y},\\[1ex]
g_1W_1 \cap g_2W_1=
(g_1 \cap g_2t^{-1})\widetilde{Y} \subseteq \widetilde{W}~~
(g_1 \neq g_2 \in [G:G_1]).
\end{array}$$
For any cell $D \subseteq \widetilde{W}$
$$U(D)=\begin{cases}
\{g_1\}&\hbox{\!\!\!if $g_1D \subseteq
\widetilde{X}_1-\bigcup\limits_{h \in [G_1;H]}(h\widetilde{Y}\cup ht\widetilde{Y})$
for some $g_1 \in [G:G_1]$}\\
\{g_1,g_2,h\}&\hbox{\!\!\!if $hD \subseteq \widetilde{Y}\times [0,1]$ for
some $h=g_1\cap g_2t^{-1} \in [G_1;H]$.}
\end{cases}$$
(ii)\qua If $W_1$ is a fundamental domain for any connected
CW complex $W$ with $\pi_1(W)=G=G_1*_H\{t\}$ then
$W=X_1\cup_{Y \times \{0,1\}}Y \times [0,1]$
as in (i)\qua, with
$$X_1=W_1/G_1,~Y=(W_1 \cap tW_1)/H.$$}
\end{Example}

\begin{Definition} {\rm Suppose that $W$ is $n$--dimensional.
Lift each cell $D^r \subseteq W$ to a cell
$\widetilde{D}^r \subseteq \widetilde{W}$.
A sequence $U=\{U_n,U_{n-1},\dots,U_1,U_0\}$ of subtrees
$U_r \subseteq T$ is {\it realized} by $W$ if the subspace
$$W(U)_1=\bigcup\limits^n_{r=0}\bigcup\limits_{D^r \subset W}
\bigcup\limits_{g_1 \in U^{(0)}_r}g_1\widetilde{D}^r
\subseteq \widetilde{W}$$
is a connected subcomplex, in which case $W(U)_1$ is a domain for
$\widetilde{W}$ with
$$W(U)_1\cap tW(U)_1=
\bigcup\limits^n_{r=0}\bigcup\limits_{D^r \subset W}
\bigcup\limits_{h \in U^{(1)}_r}h\widetilde{D}^r \subseteq \widetilde{W}$$
a connected subcomplex. Thus $U$ is realized by $C(\widetilde{W})$ and
$$C(W(U)_1)=C(\widetilde{W}(U)_1 \subseteq j_1^!C(\widetilde{W})$$
is the domain for $C(\widetilde{W})$ given by
$C_r(\widetilde{W})_1(U_r)$ in degree $r$.}
\end{Definition}

\begin{Proposition} \label{lemma2}$\phantom{9}$

{\rm (i)}\qua For any domain $W_1$ there is defined a homotopy pushout
$$\xymatrix@C-10pt@R-10pt{
(W_1\cap tW_1)/H\times \{0,1\}\ar[rr]^-{e_1} \ar[dd]_-{e_2} & &
W_1/G_1 \ar[dd]^-{f_1} \\
& \Phi & \\
(W_1\cap tW_1)/H \times [0,1] \ar[rr]^-{f_2} && W}$$
with $e_1=i_1\cup i_2\vert$, $e_2=i_3\vert$, $f_1=j_1\vert$, $f_2=j_2\vert$.
The connected 2--sided CW pair
$$(X,Y)=({\mathcal M}(e_1,e_2),(W_1 \cap tW_1)/H \times \{1/2\})$$
is a Seifert--van Kampen splitting of $W$, with a homotopy equivalence
$$f=f_1\cup f_2\co X={\mathcal M}(e_1,e_2)
\xymatrix{\ar[r]^-{\displaystyle{\simeq}}&}W.$$
{\rm (ii)}\qua The commutative square of covering projections
$$\xymatrix{ \widetilde{W}/H\times \{0,1\} \ar[r]^-{i_1\cup i_2}
\ar[d]_-{i_3} &
\widetilde{W}/G_1 \ar[d]^-{j_1} \\
\widetilde{W}/H \times [0,1] \ar[r]^-{j_2} & W}$$
is a homotopy pushout. The connected 2--sided CW pair
$$(X(\infty),Y(\infty))=({\mathcal M}(i_1\cup i_2,i_3),\widetilde{W}/H \times \{0\})$$
is a canonical injective infinite Seifert--van Kampen splitting of $W$,
with a homotopy equivalence $j=j_1\cup j_2\co X(\infty) \to W$ such that
$$\pi_1(Y(\infty))=H \subseteq \pi_1(X(\infty))=G_1*_H\{t\}.$$
{\rm (iii)}\qua For any (finite) sequence $U=\{U_n,U_{n-1},\dots,U_0\}$
of subtrees of $T$ realized by $W$ there is defined a homotopy pushout
$$\xymatrix{ Y(U)\times \{0,1\}\ar[r]^-{e_1} \ar[d]_-{e_2} &
X(U)_1 \ar[d]^-{f_1} \\
Y(U) \times [0,1] \ar[r]^-{f_2} & W}$$
with
$$\begin{array}{l}
Y(U)=(W(U)_1 \cap tW(U)_1)/H,~X(U)_1=W(U)_1/G_1,\\[1ex]
e_1=i_1 \cup i_2\vert,~e_2=i_3\vert,~f_1=j_1\vert,~f_2=j_2\vert.
\end{array}$$
Thus
$$(X(U),Y(U))=({\mathcal M}(e_1,e_2),Y(U) \times \{1/2\})$$
is a (finite) Seifert--van Kampen splitting of $W$.

{\rm (iv)}\qua The canonical infinite domain of a finite CW complex $W$
with $\pi_1(W)=G_1*_H\{t\}$ is a union of finite domains $W(U)_1$
$$\widetilde{W}=\bigcup \limits_UW(U)_1$$
with $U$ running over all the finite sequences realized by $W$.
The canonical infinite Seifert--van Kampen splitting is thus a
union of finite Seifert--van Kampen splittings
$$(X(\infty),Y(\infty))=\bigcup\limits_U(X(U),Y(U)).\eqno{\qed}$$
\end{Proposition}

\begin{Proposition}
Let $(X,Y)$ be a finite connected 2--sided CW pair
with $X=X_1\cup_{Y \times \{0,1\}}Y \times [0,1]$ for connected $X_1,Y$,
together with an isomorphism
$$\pi_1(X)=\pi_1(X_1)*_{\pi_1(Y)}\{t\}
\xymatrix{\ar[r]^-{\displaystyle{\cong}}&} G=G_1*_H\{t\}$$
preserving the HNN structures, with the structure on $G$
injective. It is possible to attach a finite number of
2-- and 3--cells to the finite Seifert--van Kampen splitting $(X,Y)$ of $X$
to obtain a finite injective Seifert--van Kampen splitting $(X',Y')$ with
$X'=X'_1\cup_{Y' \times \{0,1\}}Y' \times [0,1]$ such that
\begin{itemize}
\item[\rm (i)] $\pi_1(X')=G$, $\pi_1(X'_1)=G_1$, $\pi_1(Y')=H$,
\item[\rm (ii)] the inclusion $X \to X'$ is a homotopy equivalence,
\item[\rm (iii)] the inclusion $X_1 \to X_1'$
is a $\Z[G_1]$--coefficient homology equivalence,
\item[\rm (iv)] the inclusion $Y \to Y'$ is a $\Z[H]$--coefficient
homology equivalence.
\hfill$\qed$
\end{itemize}
\end{Proposition}

This completes the proof of the Combinatorial Transversality Theorem
for HNN extensions.

\Addresses\recd
\end{document}

%% file: gtmonout.tex
\def\ifplaintex{\expandafter\ifx\csname documentclass\endcsname\relax}

\ifplaintex 
\else
\fi

\def\gt{{\mathsurround=0pt\it $\cal G\mskip-2mu$eometry \&\ 
$\cal T\!\!$opology}}        

\def\gtm{{\mathsurround=0pt\it $\cal G\mskip-2mu$eometry \&\ 
$\cal T\!\!$opology $\cal M\mskip-1mu$onographs}}    

\def\gtp{{\mathsurround=0pt\it $\cal G\mskip-2mu$eometry \&\ 
$\cal T\!\!$opology $\cal P\!$ublications}}  

\def\recd{{\small Received:\qua\receiveddate\ifx\reviseddate\relax
\else\qquad Revised:\qua\reviseddate\fi\par}} 

\def\volumenumber#1{\def\thevolumenumber{#1}}
\def\volumeyear#1{\def\thevolumeyear{#1}}
\def\volumename#1{\def\thevolumename{#1}}
\def\papernumber#1{\def\thepapernumber{#1}}
\def\pagenumbers#1#2{\def\startpage{#1}\def\finishpage{#2}}
\def\published#1{\def\publishdate{#1}}
\def\received#1{\def\receiveddate{#1}}
\def\revised#1{\def\reviseddate{#1}}
\def\accepted#1{\def\accepteddate{#1}}
\def\asciititle#1{\def\theasciititle{#1}}

\long\def\asciiabstract#1{\long\def\theasciiabstract{#1}}

\let\\\par
\let\thevolumenumber\relax\let\thepapernumber\relax
\let\thevolumeyear\relax\let\startpage\relax
\let\finishpage\relax\let\publishdate\relax\let\receiveddate\relax
\let\reviseddate\relax\let\accepteddate\relax\let\theasciititle\relax
\let\theasciiauthors\relax
\let\theasciiabstract\relax

\let\theerratum\relax\let\theasciiemail\relax
\let\theshortauthors\relax\let\theshorttitle\relax

\def\startpage{1}\def\finishpage{15}\def\thepapernumber{77}

\volumenumber{2}
\volumename{Proceedings of the Kirbyfest}
\volumeyear{1999}

\long\def\maketitlep{   

\count0=\startpage

\gtm\nl        
{\small Volume \thevolumenumber: \thevolumename\nl 
\ifx\theerratum\relax\else Erratum \erratumnumber\nl\fi
Pages \startpage--\finishpage\nl}

\vglue 0.1truein   

{\parskip=0pt\leftskip 0pt plus 1fil\def\\{\par\smallskip}{\ifplaintex\large
\else\Large\fi\bf\thetitle}\par\medskip}   
\vglue 0.05truein 

{\parskip=0pt\leftskip 0pt plus 1fil\def\\{\par}{\sc\theauthors}
\par\medskip}%
 
\vglue 0.03truein 

{\small\leftskip 25pt\rightskip 25pt{\bf Abstract}\stdspace\theabstract

{\bf AMS Classification}\stdspace\theprimaryclass
\ifx\thesecondaryclass\relax\else; \thesecondaryclass\fi\par
{\bf Keywords}\stdspace \thekeywords\par}\vglue 7pt

}   

\font\phead=cmsl9 scaled 950
\font=cmsl9 scaled 1050
\font\pnum=cmbx10 scaled 913
\font\lnum=cmbx10 
\font\pfoot=cmsl9 scaled 950
\font=cmsl9 scaled 1050
\ifplaintex
\headline{\vbox to 0pt{\vskip -4.5mm\line{\small\phead\ifnum
\count0=\startpage ISSN 1464-8997 (on line)
1464-8989 (printed) \hfill {\pnum\folio}\else\ifodd\count0\def\\{ }%
\ifx\theshorttitle\relax\thetitle\else\theshorttitle\fi\hfill{\pnum\folio}
\else\def\\{ and }{\pnum\folio}\hfill\ifx\theshortauthors\relax\theauthors
\else\theshortauthors\fi\fi\fi}\vss}}
\footline{\vbox to 0pt{\vglue 0mm\line{\small\pfoot\ifnum\count0=\startpage
Published \publishdate:\qua\copyright\ \gtp\hfill\else
\gtm, Volume \thevolumenumber\ (\thevolumeyear)\hfill\fi}\vss
}}
\else
\makeatletter
\def\@oddhead{{\small\lhead\ifnum\count0=\startpage ISSN 1464-8997 (on line)
1464-8989 (printed) \hfill {\lnum\number\count0}\else\ifodd\count0
\def\\{ }\ifx\theshorttitle\relax \thetitle \else\theshorttitle\fi\hfill
{\lnum\number\count0}\else\def\\{ and }{\lnum\number\count0}
\hfill\ifx\theshortauthors\relax 
\theauthors\else\theshortauthors\fi\fi\fi}}\def\@evenhead{@oddhead}
\def\@oddfoot{\small\lfoot\ifnum\count0=\startpage Published \publishdate:\qua\copyright\ \gtp\hfill\else
\gtm, Volume \thevolumenumber\ (\thevolumeyear)\hfill\fi}
\def\@evenfoot{@oddfoot}
\makeatother
\fi

\let\maketitlepage\maketitlep
\let\makeshorttitle\maketitlepage
\let\maketitle\maketitlepage

\newwrite\gtoutfile
\long\gdef\makeheadfile{  
{\def\\{, }\def\s{ }
\immediate\openout\gtoutfile head.xxx
\immediate\write\gtoutfile{Proxy-for: \ifx\theasciiauthors\relax
\theauthors\else\theasciiauthors\fi\s<\ifx\theasciiemail\relax\theemail\else\theasciiemail\fi>}
\immediate\write\gtoutfile{\noexpand\\}
\immediate\write\gtoutfile{Authors: \ifx\theasciiauthors\relax
\theauthors\else\theasciiauthors\fi}
{\def\\{ }\immediate\write\gtoutfile{Title: \ifx\theasciititle\relax
\thetitle\else\theasciititle\fi}}
\immediate\write\gtoutfile{Subj-class: GT or SG, GR etc}
\immediate\write\gtoutfile{MSC-class: \theprimaryclass\ifx\thesecondaryclass\relax\else, \thesecondaryclass\fi}
\immediate\write\gtoutfile{Journal-ref: Geom. Topol. Monogr. \thevolumenumber\s
(\thevolumeyear) \startpage-\finishpage}
\immediate\write\gtoutfile{Comments: Published by Geometry and Topology Monographs at}
\immediate\write\gtoutfile{\s\s\s  http://www.maths.warwick.ac.uk/gt/GTMon\thevolumenumber/paper\thepapernumber.abs.html}
\immediate\write\gtoutfile{\noexpand\\}
\immediate\write\gtoutfile{}
\ifx\theasciiabstract\relax
\immediate\write\gtoutfile{\theabstract}\else
\immediate\write\gtoutfile{\theasciiabstract}\fi
\immediate\write\gtoutfile{}
\immediate\write\gtoutfile{\noexpand\\}
\immediate\write\gtoutfile{}
\immediate\closeout\gtoutfile}}  

\def\maketitlepage{\maketitlep\makeheadfile}
\let\makeshorttitle\maketitlepage
\let\maketitle\maketitlepage